\RequirePackage{etoolbox}
\csdef{input@path}{{style/}{graphics/}}
\documentclass[aop,MSNbibl,seceqn,dvips]{arximspdf}
\usepackage{mathrsfs,mathbh}
\usepackage{graphicx}

\doi{10.1214/13-AOP881} 
\volume{43}
\issue{2}
\pubyear{2015}
\firstpage{528}
\lastpage{571}

\makeatletter
\newcommand{\overset}{\stackrel}
\newcommand{\underset}[2]{\mathop{#2}_{#1}}
\newtheorem{theorem}{Theorem}[section]
\newproclaim{assume}[theorem]{Assumption}
\newtheorem{corollary}[theorem]{Corollary}
\newtheorem{lemma}[theorem]{Lemma}
\newproclaim{definition}[theorem]{Definition}
\newtheorem{proposition}[theorem]{Proposition}
\newproclaim{fact}[theorem]{Fact}
\newproclaim{remark}[theorem]{Remark}
\def\eqref#1{(\ref{#1})}
\newcommand{\R}{\mathbb{R}}
\newcommand{\C}{\mathbb{C}}
\newcommand{\Z}{\mathbb{Z}}
\newcommand{\N}{\mathbb{N}}

\def\SS{\mathcal{S}}
\def\floor#1{\lfloor{#1}\rfloor}

\def\eps{\varepsilon}
\def\P{\mathbb{P}}
\def\E{\mathbb{E}}

\def\md{|}
\def\Bb#1#2{{\def\md{|}#1[#2]}}

\def\Pb{\Bb\P}
\def\Eb{\Bb\E}

\def\FK#1#2#3{{\def\md{\bigm| } \mathbb{P}_{#1}^{ #2}  [  #3 ]}}
\def\EFK#1#2#3{{\def\md{\bigm| } \mathbb{E}_{#1}^{ #2}  [  #3 ]}}

\def\1{\mathbh{1}}
\def\QUAD{\mathcal{Q}} 
\def\HH{\mathscr{H}} 
\def\T{\mathcal{T}} 
\def\Sob{\mathcal{H}} 
\def\A{\mathcal{A}}
\def\<#1{\langle #1\rangle}
\def\DuminilSixArm{DuminilSixArm} 
\def\LedouxTalagrand{MR2814399}
\def\Wer09{MR2523462}
\def\Dubedat{MR2525778}
\def\SmirnovFKI{MR2680496}
\def\KSIII{KemppainenSmirnovIII}
\def\GrimmettFK{MR2243761}

\def\GPSb{GPS2b}
\def\SchrammSmirnov{SSblacknoise}
\def\SS{\SchrammSmirnov}
\def\CamiaNewmanFull{MR2249794}
\def\GHS{MR0266507}

\def\Wu{Wu}
\def\WuB{WuB}
\def\Chuck{Chuck}
\def\CamiaNewmanIsing{MR2504956}
\def\Camia{Camia}
\def\CGNproperties{CGNproperties}

\def\CGNNC{CGNNC}
\makeatother

\begin{document}
\begin{frontmatter}

\title{Planar Ising magnetization field I. Uniqueness of the
critical scaling limit}
\runtitle{Ising field I: Unique limit}

\begin{aug}
\author[A]{\fnms{Federico} \snm{Camia}\corref{}\ead[label=e1]{f.camia@vu.nl}\thanksref{T1}},
\author[B]{\fnms{Christophe} \snm{Garban}\ead[label=e2]{christophe.garban@ens-lyon.fr}\ead[label=u2,url]{http://perso.ens-lyon.fr/christophe.garban/}\thanksref{T2}}\\
\and
\author[C]{\fnms{Charles M.} \snm{Newman}\ead[label=e3]{newman@cims.nyu.edu}\ead[label=u3,url]{http://www.math.nyu.edu/faculty/newman/}\thanksref{T3}}
\runauthor{F. Camia, C. Garban and C. M. Newman}
\affiliation{VU University Amsterdam and NYU Abu Dhabi, ENS Lyon, CNRS,\\
and NYU (Courant, Shanghai) and University of California, Irvine}
\address[A]{F. Camia\\
Department of Mathematics\\
Vrije Universiteit\\
De Boelelaan 1081a\\
1081 HV Amsterdam\\
The Netherlands\\
and\\
NYU Abu Dhabi \\
\printead{e1}}
\address[B]{C. Garban\\
ENS de Lyon, CNRS\\
site Monod\\
46 all\'ee d'Italie\\
69364 LYON cedex 07\\
France\\
\printead{e2}\\
\printead{u2}}
\address[C]{C. M. Newman\\
NYU--Courant Institute of Mathematical Sciences\\
New York University\\
251 Mercer St.\\
New York, New York 10012\\
USA\\
\printead{e3}\\
\printead{u3}}

\end{aug}
\thankstext{T1}{Supported in part by NWO Grant Vidi 639.032.916.}
\thankstext{T2}{Supported in part by ANR Grant MAC2 10-BLAN-0123.}
\thankstext{T3}{Supported in part by NSF Grants OISE-0730136 and DMS-10-07524.}

\received{\smonth{10} \syear{2012}}
\revised{\smonth{9} \syear{2013}}

%
\begin{abstract}
The aim of this paper is to prove the following result. Consider the
critical Ising model on the rescaled grid $a \mathbb{Z} ^2$,
then the renormalized magnetization field
\[
\Phi^a:= a^{15/8} \sum_{x\in a \mathbb{Z}^2}
\sigma_x \delta_x,
\]
seen as a random distribution (i.e., generalized function) on the
plane, has a~unique scaling limit as the mesh size $a\searrow0$.
The limiting field is conformally covariant.
\end{abstract}

%
\end{frontmatter}

\setcounter{footnote}{3}
\section{Introduction}

\subsection{Overview}

The Ising model, introduced by Lenz in 1920 to describe ferromagnetism,
is one the most studied models of statistical mechanics. Its two-dimensional
version has played a special role in the theory of critical phenomena since
Peierls famously proved, in 1936, that it undergoes a phase transition, and
Onsager presented, in 1944, his derivation of the free energy~\cite{Onsager}.
The phase transition of the two-dimensional Ising model has been extensively
studied by both physicists and mathematicians, becoming a prototypical example
and a test case for developing ideas and techniques and for checking hypotheses.
Its analysis has helped to test
one of the fundamental beliefs of statistical mechanics that a physical system
near the critical point of a continuous phase transition is
characterized by a
single length scale, the \textit{correlation length}, which provides
the natural
length scale for the system, and that the correlation length diverges
at the critical
point. Furthermore, close to criticality, this divergence is assumed to
be solely
responsible for singularities in the thermodynamic functions; it also
implies that
the critical system has no characteristic length and is therefore
invariant under
scale transformations. This in turn suggests that all thermodynamic functions
at criticality are homogeneous functions, and predicts the appearance
of power
laws. It also means that it should be possible to rescale the critical system
appropriately and obtain a continuum model (the \emph{continuum scaling limit})
which may have more symmetries (and be therefore easier to study) than the
original discrete model, defined on a lattice. This idea is at the
heart of the
renormalization group philosophy.

Indeed, thanks to the work of Polyakov~\cite{polyakov} and others~\cite
{bpz1,bpz2},
it was understood by physicists since the early seventies that, once an
appropriate
continuum scaling limit is taken, critical statistical mechanical
models should acquire
conformal invariance, as long as the discrete models have ``enough''
rotation invariance.
This property gives important information, enabling the determination
of two- and
three-point correlation functions at criticality, when they are nonvanishing.
Because the conformal group is in general a finite dimensional Lie
group, the resulting constraints are limited in number; however,
the situation becomes particularly interesting in two dimensions,
since there every analytic function 
$f$ defines a conformal transformation, 
provided that $f'$ is nonvanishing. As a consequence, the conformal group
in two dimensions is infinite-dimensional.

After this observation was made, a large number of critical problems
in two dimensions were analyzed using conformal methods, which were applied,
among others, to Ising and Potts models, Brownian motion, the
self-avoiding walk,
percolation and diffusion limited aggregation.
The large body of knowledge and \mbox{techniques} that resulted, starting
with the work of Belavin, Polyakov and Zamolodchikov~\cite{bpz1,bpz2}
in the early eighties, goes under the name of Conformal Field Theory (CFT).
In two dimensions, one of the main goals of CFT and its most
important application to statistical mechanics is a complete
classification of all universality classes via irreducible representations
of the infinite-dimensional Virasoro algebra (see, e.g., \cite{DiMS97}).

CFT has proved very powerful, but it also has limitations.
First of all, the theory deals primarily with correlation functions
of \textit{local} (or quasi-local) operators, and is therefore not always
the best tool to investigate global quantities.
Secondly, given some critical lattice model, there is no way, within
the theory itself, of deciding to which CFT it corresponds.
A third limitation, 
at least from a mathematician's perspective, is its lack of
mathematical rigor.

Quite remarkably, some of the most recent and significant developments
in the area of
two-dimensional critical phenomena 
have emerged in the mathematics literature, using new mathematical
tools that are free
from at least some of the limitations of CFT. These tools have
permitted to rigorously establish
the conformal invariance of several models and prove various
results/conjectures that had
first appeared in the physics literature, as well as novel results that
have shed new light on
the theory of two-dimensional critical phenomena.

In 1999, Aizenman and Burchard \cite{ab}, based on earlier work of
Aizenman \mbox{\cite{Aizenman95,Aizenman98}},
proposed a framework for proving tightness, and thus the existence of
subsequential scaling limits for the
distribution of random paths in the scaling limit. Their results found
applications in much of the subsequent work
on scaling limits of interfaces.

In 2000 and 2001, Kenyon \cite{Kenyon00,Kenyon01} proved conformal
invariance of the two-dimensional
dimer lattice model (or domino tiling model) in the scaling limit, and
related the latter to the Gaussian free field.

A very significant breakthrough was the introduction by Schramm \cite
{SchrammSLE}
of the Stochastic (Schramm--)Loewner Evolution (SLE) and its subsequent analysis
and application to the scaling limit problem for several models, most
notably by
Lawler, Schramm and Werner \cite{lsw04}, and by Smirnov \cite
{Smirnov01} (see also \cite{cn07}).
The subsequent introduction of the Conformal Loop Ensembles (CLEs)
\cite{Werner03,CN04,MR2249794,SheffieldCLE,SW12}, which are collections
of SLE-type,
closed curves, provided an additional tool to analyze the scaling limit
geometry of critical models.

Substantial progress in the rigorous analysis of the two-dimensional
Ising model at criticality
was made by Smirnov \cite{MR2680496} with the introduction and scaling
limit analysis of the
``fermionic observables,'' also known as ``discrete holomorphic
observables'' or
``holomorphic fermions.'' (Similar objects had been considered by
Mercat \cite{Mercat01}
and had appeared in the physics literature---see \cite{KC,CR}.)
These have proved extremely useful in studying the Ising model in
finite geometries with
boundary conditions and in establishing the conformal invariance of the
scaling limit of
various quantities, including the energy density \cite
{hs10,hongler-thesis} and spin
correlation functions \cite{arXiv1202.2838}.
(An independent derivation of the critical Ising correlation functions
in the plane was
obtained in \cite{arXiv1112.4399}.)

The result of Chelkak, Hongler and Izyurov \cite{arXiv1202.2838} on the
scaling limit of
spin correlation functions is the main ingredient in our second proof
of the uniqueness
of the scaling limit of the Ising magnetization, presented in
Section~\ref{second}, and
it is also used in our first proof, presented in Section~\ref{first}.

Our second proof essentially consists in showing the existence of the
scaling limit of the
characteristic function of the discrete field.
Our first derivation is very different in spirit from the second; it is
more geometric in nature
and is based on the RSW-type result for FK-Ising percolation of
Duminil-Copin, Hongler and
Nolin \cite{arXiv0912.4253}, and on scaling limit results for FK-Ising
percolation
\cite{Kemppainen-thesis,CDHKS,KS1,KemppainenSmirnovIII}. This is in fact a
conditional
proof of uniqueness since it relies on a scaling limit result that,
although very plausible,
does not follow immediately from known results (see Section~\ref{full-s-lim} for a
detailed explanation).

Earlier influential results on the scaling limit of the two-dimensional
Ising model that
have been a source of inspiration for the present paper include those
of Abraham
\cite{Abraham73,Abraham78}, Abraham and Reed \cite{AR76}, De Coninck
\cite{DC87},
De Coninck and Newman \cite{DCN90} and Palmer \cite{Palmer-book}.

\subsection{Definitions and main results}\label{s.renormalization}

Consider the Ising model on the\break rescaled grid $a \Z^2$ at the \textit
{critical temperature} $\beta_c=\beta_c(\Z^2)$,
with zero external magnetic field.
(We refer, e.g., to \cite{\GrimmettFK} for a nice introduction to
the Ising model.) We will be interested in the following object.

\begin{definition}\label{def:magnetization-field}
The renormalized \textit{magnetization field} $\Phi^a$, a random
distribution on the plane, is
\[
\Phi^a:= \Theta_a \sum_{x\in a \Z^2}
\sigma_x \delta_x,
\]
where $\Theta_a$ is a well-chosen renormalization factor.
(In fact, we will use a slightly modified version; see Definition~\ref{d.MA2}.)

In most of the rest of the paper,
we will fix\footnote{This particular choice assumes Wu's result \cite{\Wu,\WuB}. Note that
this choice may be debatable. For example,
the authors of \cite{arXiv1202.2838} do not assume Wu's result. Without such
an assumption, our results remain valid with
$\Theta_a$ defined more implicitly.
See Remark~\ref{remark-wu} and Section~\ref{Theta}. 
}
$\Theta_a:= a^{15/8}$.
\end{definition}

In the scaling limit, the magnetization field is expected to converge
to a Euclidean random field corresponding to the simplest
reflection-positive conformal field theory \cite{bpz1,bpz2} (see also
\cite{SimonPphi2}).
Our main theorem can be stated as follows.

%
\begin{theorem}[(Scaling limit)]\label{th.main}
The magnetization field $\Phi^a$ converges in law as the mesh size
$a\searrow0$ to a limiting random distribution $\Phi^\infty$.
The convergence in law holds in the Sobolev space $\Sob^{-3}$ under the
topology given by $\|\cdot\|_{\Sob^{-3}}$.
(See Appendix~\ref{tight}.)
\end{theorem}

In fact, our result holds for any bounded simply connected domain
$\Omega\subset\C$ with a smooth boundary.\footnote{In principle,
the results do not require the domain to be simply connected and have a
smooth boundary, but these assumptions allow us to directly
use various results from the literature that so far were proved only in
the simply connected case with smooth boundary.}
More precisely, consider a simply connected domain $\Omega$ in the
plane which contains the origin, and let $\Omega_a$ denote its
approximation by the grid $a \Z^2$ of mesh size $a$, that is, $\Omega
_a:= \Omega\cap a \Z^2$. (The approximation might not be simply
connected anymore, so in this case we keep only the connected component
of the origin.) Consider the Ising model in $\Omega_a$ with $+$ boundary
conditions, at the critical temperature $\beta_c$ (we will also analyze
the case of free boundary conditions). The above definition of
(renormalized) magnetization field easily extends to this setting:
\[
\Phi^a_{\Omega} \bigl(= \Phi^a \bigr):=
\Theta_a \sum_{x\in\Omega_a} \sigma_x
\delta_x.
\]



\begin{theorem}\label{th.main2}
Let $\Omega$ be a bounded simply connected domain of the plane with a
smooth boundary. Consider the critical Ising model with $+$ or free
boundary conditions in $\Omega_a$. Then the magnetization field $\Phi
^a_\Omega=\Phi^a$ converges in law as the mesh size $a\searrow0$ to
a limiting random distribution $\Phi^\infty_\Omega= \Phi^\infty$.
The convergence in law holds in the Sobolev space $\Sob^{-3}=\Sob
^{-3}(\Omega)$ under the topology given by
$\|\cdot\|_{\Sob^{-3}}$. (See Appendix~\ref{tight}.)
\end{theorem}

We now explain how to choose the scaling factor $\Theta_a$.
For any bounded domain~$\Omega$, the magnetization $M^a:=\sum_{x\in
\Omega_a} \sigma_x$,
where $\sigma_x$ denotes the Ising \textit{spin variable} at $x$, has variance
\begin{eqnarray*}
\operatorname{Var} \bigl[M^a \bigr] & =& \sum
_{x,y\in\Omega_a^2} \Bb\E{\sigma_x \sigma_y}.
\end{eqnarray*}
It can be shown (see Proposition~\ref{pr.1and2} in Appendix~\ref{s.appendix}) that\footnote{In this paper, $f(a) \asymp g(a)$ as $a\searrow0$ means that
$f(a)/g(a)$ is bounded away from
$0$ and $\infty$ while $f(a) \sim g(a) $ means that $f(a)/g(a) \to1$
as $a\searrow0$.} as $a\searrow0$,
\[
\operatorname{Var} \bigl[M^a \bigr] \asymp a^{-4}
\mathbb{E}_{\mathbb{C}_{a}}[ \sigma_{0_a}\sigma_{(\sqrt{2} + \sqrt{2}i)_a}],
\]
where $\mathbb{C}_{a}$ denotes the square grid $a\mathbb{Z}^2$ and
$0_a$ and $(\sqrt{2} + \sqrt{2}i)_a$ stand for
lattice approximations of the points $0,\sqrt{2} + \sqrt{2}i \in\mathbb{C}$.

Following the notation of \cite{arXiv1202.2838}, let us introduce the quantity
%
%
\begin{equation}
\label{e.varrho} \varrho(a ):=\mathbb{E}_{\mathbb{C}_{a}}[\sigma_{0_a}
\sigma_{(\sqrt{2} + \sqrt{2}i)_a}].
\end{equation}
From the above discussion, it is thus natural to scale our
magnetization field by a scaling factor of order $a^{2} \varrho(a
)^{-1/2}$. In most of the rest of this paper (until Section~\ref{Theta}), we will assume the following celebrated result by Wu.

%
\begin{theorem}[(Wu, see \cite{\WuB,\Wu})]\label{th.Wu}
There exists an explicit constant $c>0$ such that as $a\searrow0$
%
%
\begin{equation}
  \varrho(a ) \sim c a^{1/4}.
\end{equation}
\end{theorem}

Assuming this asymptotic result leads to the choice
%
%
\begin{equation}
  \Theta_a:= a^{15/8}
\end{equation}
for the scaling factor of the magnetization field.


%
\begin{remark}\label{remark-wu}
We believe that it is reasonable to assume Wu's result since it is
considered to be among the rigorous results obtained in the theoretical
physics literature (yet, according to some experts, although there is
no theoretical gap, some details need to be filled in). Nevertheless,
this choice may be debatable. In \cite{arXiv1202.2838}, for example, the
authors decided to state their result without assuming Wu's result. In
our case, if one avoids assuming Wu's asymptotic, all our results
remain valid by replacing the above
formula for $\Theta_a$ by the more cumbersome $\Theta_a:= a^{2} \varrho
(a )^{-1/2}$; see Section~\ref{Theta}.
\end{remark}

The following corollary on the renormalized magnetization random
variable follows easily from Theorem~\ref{th.main2}
when taking $\Theta_a:= a^{15/8}$.
We state it in a similar way as one would state a central limit
theorem. Indeed, dealing with a sum
of random variables, the result below has a classical flavor; its
interest lies in the strong dependence of the random
variables being added which leads to a nondegenerate non-Gaussian
limit. (Note that the choices of the unit square
as domain and of $+$ boundary conditions are made only for concreteness
and are not essential.)

%
\begin{corollary}\label{c.CLT}
Consider the critical Ising model in the $N\times N$ square $\Lambda_N$
with $+$ boundary conditions on $\partial\Lambda_N$.
Then the random variable
\[
\frac{1} {N^{15/8}} \sum_{x\in\Lambda_N} \sigma_x
\]
converges in law as $N\to\infty$.
\end{corollary}

The result in Corollary~\ref{c.CLT} can be expressed also in terms of
the renormalized \textit{magnetization},
defined below. The magnetization is the \textit{order parameter} of the
Ising phase transition, that is the extra parameter
of the model that is needed, due to the spontaneous breaking of the
spin symmetry below the critical temperature, to describe
the thermodynamics of the low-temperature phase. A fundamental belief
of statistical mechanics is that, near the phase
transition point, the order parameter is the only important
thermodynamic quantity.

%
\begin{definition}[(Renormalized magnetization)]\label{d.AM}
For any simply connected domain $\Omega$ with boundary condition $\xi$
on $\partial\Omega$
(which in this paper will be either $+$ or $\mathrm{free}$), let
$m_\Omega^a$ be the renormalized magnetization in the domain
$\Omega$ defined by
\begin{eqnarray*}
  m_\Omega^a &:=& \bigl\langle
\Phi^a, 1_\Omega\bigr\rangle
\\
&=& \Theta_a \sum_{x\in\Omega_a}
\sigma_x
\\
&=& a^{15/8} \sum_{x\in\Omega_a}
\sigma_x.
\end{eqnarray*}
\end{definition}

Exactly as in Corollary~\ref{c.CLT}, $m^a_\Omega$ converges in law to a
limit $m^\infty_{\Omega}$ as $a\searrow0$.
(The limiting law will depend on the boundary condition $\xi$.)

The limiting magnetization field is non-Gaussian: this can be seen from
the correlation functions computed by
Chelkak, Hongler and Izyurov \cite{arXiv1202.2838}, which do not
satisfy Wick's formula.
We will list a number of other properties satisfied by the limiting
fields $\Phi^\infty_\C$, $\Phi^\infty_\Omega$, $m^\infty_{\Omega}$
in Section~\ref{properties}. Most of them will be proved in \cite
{\CGNproperties}. In particular, we will prove in \cite{\CGNproperties}
that the tail behavior of these limiting fields is of the form $\exp(-
c x^{16})$.

We conclude this section by stating the conformal covariance properties
of the scaling limit $\Phi^\infty$
of the lattice magnetization field.


%
\begin{theorem}[(Conformal covariance of $\Phi^\infty$)]\label{th.CC}
Let $\Omega, \widetilde\Omega$ be two simply connected domains of the
plane (not equal to $\C$) and let $\phi\dvtx \Omega\to\widetilde\Omega$
be a conformal map. Let
$\psi=\phi^{-1}$ be the inverse conformal map from $\widetilde\Omega
\to\Omega$.
Let $\Phi^\infty$ and $\widetilde\Phi^\infty$ be the continuum
magnetization fields, respectively, in $\Omega$, $\widetilde\Omega$.
Then the pushforward distribution $\phi*\Phi^\infty$ of the\vspace*{1pt} random
distribution $\Phi^\infty$ has the same law as the random distribution
$|\psi'|^{15/8} \widetilde\Phi^\infty$, where the latter distribution
is defined as
\[
\bigl\langle \bigl|\psi'\bigr|^{15/8} \widetilde\Phi^\infty,
\tilde f\bigr\rangle:= \bigl\langle \widetilde\Phi^\infty, w\mapsto\bigl|
\psi'\bigr|^{15/8}(w) \tilde f(w)\bigr\rangle
\]
for any test function $\tilde f\dvtx \widetilde\Omega\to\C$.
\end{theorem}

In the particular case of the renormalized magnetization in squares of
various scales, the above conformal covariance property
can be expressed as follows.

%
\begin{corollary}\label{c.SI}
Let $m^\infty$ be the scaling limit of the renormalized magnetization
in the square (i.e., $m^\infty= \<{\Phi^\infty, 1_{[0,1]^2}}$).
For any $\lambda>0$, let $m^\infty_\lambda$ be the scaling limit of the
renormalized magnetization in the square $[0,\lambda]^2$. Then one has
the following identity in law:
%
%
\begin{equation}
  m^\infty_\lambda\overset{d} {=}
\lambda^{15/8} m^\infty.
\end{equation}
\end{corollary}

\subsection{Brief outline of the proofs}
We will give two proofs of our main result, Theorem~\ref{th.main2},
each proof having its own advantages.
Let us briefly sketch in this subsection our two strategies. They both
start with the same tightness step.

\subsubsection{Tightness}
Tightness of the random variable $m^a_{\Omega}:= \<{\Phi^a, 1_{\Omega
}}$ was already proved in \cite{\CamiaNewmanIsing}
(see also \cite{\Camia}).
To prove tightness for the random distribution~$\Phi^a$ we choose to
work on the Sobolev space $\Sob^{-3}$ and use the setup of
\cite{MR2525778}.
With this setup the proof is relatively standard but somewhat
technical, so we give an outline here and present the details in
Appendix \ref{tight} at the end of the paper.
Below and in the rest of the paper, except for Appendix \ref{tight}, we restrict
our attention to the magnetization in the unit square $[0,1]^2$, in
order to simplify the notation.
The extensions to other domains and to the full plane are discussed in
Appendix \ref{tight}.


For any $a>0$, we will consider our magnetization field $\Phi^a$ as an
element of the Polish space $\Sob^{-3}$ with operator norm $\| \cdot\|
_{\Sob^{-3}}$
(see Appendix \ref{tight} for precise definitions).
Since Dirac point masses do not belong to $\Sob^{-\alpha}$ for $\alpha
\leq1/2$, it will be convenient to change slightly the definition
of the distribution $\Phi^a$ to the following definition.

\begin{definition}\label{d.MA2}
We let
\[
\Phi^a:= a^{15/8} \sum_{x\in[0,1]^2\cap a\Z^2}
\frac{\sigma_x} {a^2} 1_{S_a(x)},
\]
where $S_a(x)$ denotes the square centered at $x$ of side-length $a$.
\end{definition}

With this new definition, $\Phi^a$ belongs to $L^2$, and hence has a
Fourier expansion. Using the latter, it is not hard to show that
\[
\limsup_{a\searrow0}
\E\bigl[ \bigl\| \Phi^a \bigr\|_{\Sob^{-2}}^2 \bigr] < \infty,
\]
provided that the boundary condition on the square $[0,1]^2$ consists
of finitely many arcs of $+,-$ or \textit{free} type.
This is enough to prove that $(\Phi^a)_{a>0}$ is \textit{tight} in the
space $\Sob^{-3}$ thanks to the Rellich theorem,
which implies that, for any $R>0$, the ball
\[
\overline{B_{\mathcal{H}^{-2}}(0,R)}
\]
is compact in $\mathcal{H}^{-3}$.
As a consequence, we have the following proposition.

\begin{corollary}\label{cor.sub}
Consider the magnetization in the unit square $[0,1]^2$ with boundary
condition consisting of finitely many arcs of $+,-$ or \textit{free} type.
Then there is a subsequential scaling limit $\Phi^*$, that is, a random
distribution $\Phi^*\in\Sob^{-3}$ such that for a certain subsequence
$a_k \searrow0$, $\Phi^{a_k}$ converges in law to $\Phi^*$ for the
topology on $\Sob^{-3}$ induced by $\|\cdot\|_{\Sob^{-3}}$.
\end{corollary}

\subsubsection{First proof}

In the first proof (Section~\ref{first}), we rely on the FK
representation of the Ising model which allows us to decompose
the distribution $\Phi^a$ as a sum over the FK clusters, where each
cluster $C$ carries an independent random sign
$\sigma_{C}\in\{-1,1\}$.
Two important ingredients in this proof are the RSW-type result for
FK-Ising percolation of Duminil-Copin, Hongler
and Nolin \cite{arXiv0912.4253} and the $k=0$ case of Theorem 1.3 of
\cite{arXiv1202.2838}. (We note that
the use of the latter result could probably be avoided by relying on an
argument similar to the one used in \cite{arXiv1008.1378}
to prove the rotational invariance of the percolation two-point function.)
The drawback of this approach is that we need to rely on the uniqueness
of the full scaling limit of FK percolation
(see Assumption~\ref{a.uniqueness}). Note that the main argument, which
consists in constructing area measures
on critical FK clusters, is somewhat close to the construction of
``pivotal measures'' in \cite{arXiv1008.1378}.

\subsubsection{Second proof}

Our second proof (Section~\ref{second}), as opposed to the first one,
does not rely on any assumption
(besides assuming Wu's result if one wants to keep the scaling $\Theta
_a=a^{15/8}$).
For any bounded domain $\Omega$, the idea is to characterize the limit
of $\Phi^a$ by showing that the quantities
\[
  \phi_{\Phi^a}(f) = \E\bigl[e^{i\<{\Phi^a,f}}\bigr],
\]
converge as $a\searrow0$ for any test function $f\in\Sob^3$. The main
ingredients are the breakthrough results by
Chelkak, Hongler and Izyurov on the convergence of the $k$-point
correlation functions as well as our Propositions~\ref{pr.expomoments}
and~\ref{pr.close}.

\subsubsection{Proof of the conformal covariance properties}
We briefly discuss how to prove Theorem~\ref{th.CC}.
\begin{enumerate}[1.]
\item[1.] If one wants to follow the setup of our first proof (Section~\ref{first}), then the conformal covariance property is proved exactly in
the same fashion as Theorem~6.1 in~\cite{arXiv1008.1378} on the conformal
covariance of the \textit{pivotal measures} for critical percolation on
the triangular lattice, except that here
one would have conformal covariance of the ensemble of FK area measures.
\item[2.] If one wants to follow the setup of our second proof (Section~\ref{second}), then Theorem~\ref{th.CC} is even easier to obtain, since it
follows easily from the conformal covariance properties of the
$k$-point functions established in the main result, Theorem~1.3, of
\cite{arXiv1202.2838}.
\end{enumerate}

In the rest of the paper, in order to simplify the notation, we will
stick to the magnetization in the unit square $[0,1]^2$.
The extension to other domains as well as to the full plane can be done
using the methods presented in Appendix \ref{tight}.\vadjust{\goodbreak}

\section{First proof of the scaling limit of \texorpdfstring{$\Phi^a$}{Phia} using area measures on FK clusters}\label{first}

\subsection{The general strategy} \label{sec:strategy}

\begin{figure}

\includegraphics{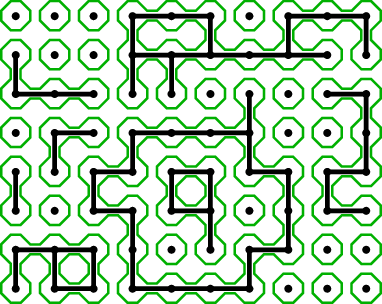}

\caption{Example of an FK bond configuration in a rectangular region.
Black dots represent vertices of $a {\mathbb Z}^2$, black horizontal
and vertical
edges represent FK bonds. The FK clusters are highlighted by lighter (green)
loops on the medial lattice.}\label{fig:loops}
\end{figure}

The FK representation of the Ising model with zero external magnetic
field is
based on the $q=2$ random-cluster measure $P_p$ (see~\cite{KF,FK,MR2243761}
for more on the random-cluster model and its connection to the Ising model).
A spin configuration distributed according to the unique
infinite-volume Gibbs
distribution with zero external magnetic field and inverse temperature
$\beta\leq\beta_c$ can be obtained in the following way. Take a random-cluster
(FK) bond configuration on $a {\mathbb Z}^2$ distributed according to
$P_p$ with
$p=p(\beta)= 1 - e^{-2\beta}$, and let $\{ C^a_i \}$ 
denote the corresponding collection of FK clusters, where a cluster
is a maximal set of vertices of the square lattice connected via bonds
of the FK configuration (see Figure~\ref{fig:loops}).
One may regard the index $i$ as taking values in the natural numbers,
but it is better to think of it as a dummy countable index without
any prescribed ordering, like one has for a Poisson point process.
Let $\{ \sigma_i \}$ 
be \mbox{($\pm1$)-}valued, i.i.d., symmetric random variables, and assign
$\sigma_x=\sigma_i$
for all $x \in C^a_i$; then the collection $\{ \sigma_x \}_{x \in
{a\mathbb Z}^2}$
of spin variables is distributed according to the unique infinite
volume Gibbs distribution
with zero external magnetic field and inverse temperature $\beta$.

%


Using the FK representation, we can write the renormalized
magnetization field~$\Phi^a$ from Definition~\ref{def:magnetization-field} as follows:
%
%
\begin{equation}
\label{eq:lattice-field} \Phi^a \stackrel{\mathrm{dist.}} {=} \sum
_i \sigma_i \mu_{C^a_i},
\end{equation}
where $\mu_{C^a_i}:= \Theta_a \sum_{x \in C^a_i}\delta_x$ and the
$\sigma_i$'s, as before, are $(\pm1)$-valued, symmetric random variables
independent of each other and everything else. We call the rescaled counting
measure $\mu_{C^a_i}$ the \textit{area measure} of the cluster $C^a_i$.

Roughly speaking, our first proof of uniqueness consists in showing
that area
measures have a scaling limit which is measurable with respect to the scaling
limit of the collection of all macroscopic crossing events. (By
crossing event we
mean the occurrence of a path of FK bonds crossing a certain domain between
two disjoint arcs on its boundary.)

To understand why this should be sufficient, let us associate in a
unique way
to each area measure $\mu_{C^a_i}$ the interface $\gamma^a_i$
in the medial lattice between the corresponding (rescaled) FK cluster $C^a_i$
and the surrounding FK clusters. 
Such interfaces form closed curves, or loops, which separate the corresponding
clusters $C^a_i$ from infinity (see Figure~\ref{fig:loops}). Announced results
for FK percolation \cite{\KSIII} identify the scaling limit of those
loops with
$\mathrm{CLE}_{16/3}$, a random collection of nested loops which are locally distributed
like $\mathrm{SLE}_{16/3}$ curves. The uniqueness of the collection of all macroscopic
crossing events should then be a consequence of the results announced
in \cite{\KSIII}
(see Assumption~\ref{a.uniqueness} below and the discussion following
it for more information).

The area measure of a cluster counts the number of vertices in that cluster.
In particular, the sum of the area measures of clusters of diameter
larger than
some $\varepsilon$ counts the number of vertices from which a path of
FK bonds
of diameter larger than $\varepsilon$ originates. We call the occurrence
of a path
of FK bonds of diameter larger than $\varepsilon$ in the scaling limit a
\textit{macroscopic
one-arm event} (see Section~\ref{ss.measurablevents} for a precise
definition of
arm events), and we will sometimes call a \textit{one-arm vertex} a vertex
from which
such a path originates. Since the area measures of macroscopic clusters count
one-arm vertices, it is reasonable to expect that they be measurable
with respect
to the collection of all macroscopic crossing events. Indeed, the
analogous result
for Bernoulli percolation was proved in \cite{arXiv1008.1378}.
Moreover, it is shown in \cite{MR2504956} that, in the scaling limit,
only area
measures corresponding to macroscopic clusters contribute to the magnetization,
therefore the latter should also be measurable with respect to the
collection of all
macroscopic crossing events.

We end this section by briefly explaining the idea behind the proof
that area measures are
measurable with respect to the collection of macroscopic crossing
events. Our proof of
this fact will follow closely the proof of \cite{arXiv1008.1378} for
Bernoulli percolation. This
can be done because the main tools used in \cite{arXiv1008.1378}, such
as FKG, RSW and
certain bounds on the probability of arm events, are also available for
Ising-FK percolation.

Although the proof is rather technical, the underlying idea is simple.
Suppose we are interested
in the Ising model on $a {\mathbb Z}^2 \cap[0,1]^2$. We superimpose on
$[0,1]^2$ a~square
grid with mesh $\eps$, with $\eps$ much smaller than 1 but much larger
than $a$.
We will show that the sum of the area measures of macroscopic FK
clusters is well approximated
by the number of squares of the $\eps$-grid that intersect a
macroscopic FK cluster, times the
``mean'' number of one-arm vertices inside a square of the \mbox{$\eps$-}grid
(we ignore all boundary
issues in this discussion). This means that the number of one-arm
vertices can be estimated by
looking only at the macroscopic features of the FK configuration (i.e.,
at the collection of macroscopic
crossing events).

%
\begin{remark}
There are several advantages to the approach presented in this section.
First, it shows that $\Phi^\infty$ is
measurable with respect to the full scaling limit of FK-Ising
percolation (plus a collection of random signs).
Furthermore, it gives us a good way to visualize the magnetization in
terms of area measures, as in \eqref{eq:lattice-field}.

Such a geometric representation as \eqref{eq:lattice-field} would also
be possible directly for the limiting
magnetization field $\Phi^{\infty}$ if one could obtain the scaling
limit of the collection of all individual,
macroscopic area measures. This should be possible with methods similar
to those used in this paper and
in \cite{arXiv1008.1378}, and the resulting scaling limit should be expressible
as a collection of orthogonal measures
supported on ``continuum FK clusters.''

We do not pursue this here since looking at the total number of $+$ and
$-$ macroscopic one-arm vertices
in a given domain is sufficient to prove the uniqueness of the limiting
magnetization field, but we point out that
approximating $\Phi^\infty$ using an a.s. finite number of signed
measures could be useful if one wanted to
determine the smallest $\eps>0$ such that $\Phi^\infty$ is in the
Sobolev space $\Sob^{-\eps}$. The latter
problem is briefly discussed in Remark~\ref{r.baralpha}.
\end{remark}

\subsection{Setup for the proof of convergence}


\subsubsection{Notation, space of percolation configurations,
compactness}\label{ss.notations}

We\break will work with the following setup: denote by $\sigma_a$ an Ising
configuration on $a \Z^2 \cap[0,1]^2$.
As explained in Section~\ref{sec:strategy}, $\sigma_a$ can be obtained
from an FK configuration
$\omega_a$ on $a \Z^2 \cap[0,1]^2$
by flipping an independent $\{\pm\}$ fair coin for each cluster of
$\omega_a$.
Let $\omega_a^+$ (resp., $\omega_a^-$) be the configuration consisting
of the clusters of $\omega_a$ which have been chosen
to be plus (resp., minus). Let us denote by $\bar\omega_a$ the coupled
pair $(\omega_a^+, \omega_a^-)$.
Note that one has $\omega_a = \omega_a^+ \cup\omega_a^-$.

It will be very convenient to consider these FK configurations $\omega
_a$ as (random) variables
in a compact metrizable space $(\HH,\T)= (\HH_{[0,1]^2}, \T)$ which
encodes all macroscopic crossing events.
We say that an FK configuration $\omega_a$ contains a crossing of a
domain $D$ between two disjoint arcs, $I_1$
and $I_2$, of its boundary $\partial D$ if there is a collection of
edges from $\omega_a$ such that the edges form a connected set
contained in $D$ except for two edges which intersect $I_1$ and $I_2$,
respectively.

The compact space $(\HH,\T)$ is not specific to our study of FK
percolation and one can in fact rely here on the
setup which was introduced by Schramm and Smirnov in \cite{\SS} in the
case of independent percolation ($q=1$).
Very briefly, it works as follows: the space of percolation
configurations built in \cite{\SS} is the space of \textit{closed
hereditary subsets} of the space of \textit{quads} $(\QUAD,d_\QUAD)$.
Roughly speaking, this means that a point $\omega\in\HH$ corresponds
to a family of quads $Q\in\QUAD$ which is closed in $(\QUAD,d_\QUAD)$
and which satisfies the following constraint: if $Q\in\omega$ and $Q'$
is ``easier'' to be traversed, then $Q'$ is in $\omega$ as well. In
\cite{\SS}, it is proved that this space $\HH$ can be endowed with a
topology $\T$ so that the topological space $(\HH, \T)$ is compact,
Hausdorff and metrizable. For convenience, we will choose a
(nonexplicit) metric $d_\HH$ on $\HH$ which induces the topology $\T$.
See \cite{\SS} for a clear exposition of the topological space $(\HH,\T
)$. See also \cite{arXiv1008.1378,\GPSb}.

Since we will need the crossing properties of the $+$ versus the $-$
clusters, we will in fact consider $\bar\omega_a=(\omega_a^+,\omega_a^-)$
as a random variable in the compact metrizable space $\HH\times\HH$
endowed with the product topology.

What is known about the limit as $a\to0$ of the coupling $\bar\omega
_a = (\omega_a^+, \omega_a^-) \in\HH\times\HH$?
First of all, the tightness for $(\bar\omega_a)_{a>0}$ follows
immediately from the compactness of $(\HH\times\HH,\T\otimes\T)$.

%
\begin{fact}
The random variable $\bar\omega_a=(\omega_a^+,\omega_a^-)$ is in $\HH
\times\HH$ (with the product topology). Since $\HH\times\HH$ is
compact for the product topology $\T\otimes\T$, there are
subsequential scaling limits for $(\bar\omega_a)_{a>0}$ as $a\to0$.
\end{fact}


Our goal will be to show that the scaling limit $\Phi^{\infty}$ of the
magnetization field $\Phi^a$
is measurable with respect to the scaling limit of $\bar\omega_a$.
This is the content of the main result of this section,
Theorem~\ref{th.uniqField}. That theorem, together with Assumption~\ref{a.uniqueness} below, immediately implies the
uniqueness of $\Phi^{\infty}$.

\subsubsection{Scaling limit for 
\texorpdfstring{$(\bar\omega_a)_{a>0}$}{(baromegaa)a>0}} \label{full-s-lim}


It is known since the breakthrough paper \cite{\SmirnovFKI} that
certain discrete ``observables'' for critical FK-percolation are
asymptotically \textit{conformally invariant}. These observables can then
be used \cite{CDHKS} to prove that interfaces have
a scaling limit described by $\mathrm{SLE}_{16/3}$ curves. In our case, we need
a \textit{full scaling limit} result. Indeed, our later results
in this section of the paper are based on the following hypothesis.

\begin{assume}\label{a.uniqueness}
The coupled configurations $\bar\omega_a=(\omega_a^+,\omega_a^-)$
considered as random points in $(\HH\times\HH, \T\otimes\T)$ have a
(unique) scaling limit as $a\searrow0$;
they converge in law to a continuum $\mathrm{FK}$ $\bar\omega_\infty=
(\omega_\infty^+,\omega_\infty^-)$.
\end{assume}

This assumption is very reasonable, based on the convergence of
discrete interfaces to $\mathrm{SLE}_{16/3}$ curves \cite{CDHKS}.
An even clearer evidence is provided by the work in progress \cite
{\KSIII},
where it is shown that the branching exploration tree converges to the
branching $\mathrm{SLE}_{16/3}$ tree.
However, as explained in \cite{\SS}, it is not always easy to go from
one notion of scaling limit to another. In the case
of Bernoulli percolation (i.e., the random-cluster model with $q=1$),
the first and third author proved \cite{\CamiaNewmanFull}
the existence and several properties of the full scaling limit as the
collection of all cluster boundaries,
building the limit object from $\mathrm{SLE}_6$ loops, and, as explained in
\cite{arXiv1008.1378}, Section~2.3,
their results imply convergence also in the ``quad topology''~$(\HH,\T)$.

In our present case, the FK percolation analog of the result contained
in \cite{\CamiaNewmanFull}, that is,
the convergence of $\omega_a$ to $\mathrm{CLE}_{16/3}$, was announced in \cite
{\KSIII}.
From this convergence result, following \cite{arXiv1008.1378},
Section~2.3, and using Corollary 5.9 of \cite
{\DuminilSixArm}
instead of the analogous result for Bernoulli percolation, one should
be able to obtain the convergence of
$\omega_a$ to $\omega_\infty$ in the topological space $(\HH,\T)$,
exactly as in the case of Bernoulli percolation.

This step would justify the convergence of $\omega_a$ to $\omega_\infty
$ in $(\HH,\mathcal{T})$, but we need slightly more,
that is, the convergence of $\bar\omega_a$ to $\bar\omega_\infty$.
However, note that the configurations $\omega_a^+$
and $\omega_a^-$ can be obtained from $\omega_a$ by tossing a fair coin
for each cluster in $\omega_a$ to decide its
sign. This suggests that the convergence of the configurations $\omega
_a^+$ and $\omega_a^-$ should follow from the
same arguments giving the convergence of $\omega_a$.

While the preceding discussion is clearly not a complete proof, it
explains why Assumption~\ref{a.uniqueness} is very
reasonable in the light of the announced results on the full scaling
limit of FK percolation.



\subsubsection{Measurable events in 
\texorpdfstring{$\mathscr{H}\times\mathscr{H}$}
{HtimesH}}\label{ss.measurablevents}

In this subsection, we follow very closely Section~2.4. in \cite{arXiv1008.1378}. We refer to that paper for more details and will only highlight
briefly how to adapt the definitions to our present case.

Let $A=(\partial_1 A, \partial_2 A) \subset[0,1]^2$ be a fixed
topological annulus whose boundary, $\partial_1 \cup\partial_2$
is composed of piecewise smooth cirves. 
We will often rely on the \textit{one-arm events} $\A^{\pm}=\A_1^\pm$
which are in the Borel sigma field of
$(\HH\times\HH, \T\otimes\T)$ and which are defined as follows:
%
%
\begin{eqnarray}
\A_1^+&:=& \bigl\{ \bar\omega\in\HH\times\HH, \exists Q
\in\QUAD\mbox{ s.t. } Q\in\omega^+
\nonumber\\[-8pt]\\[-8pt]
&&\hspace*{22pt} \mbox{ and $Q$ connects } \partial_1
A \mbox{ with } \partial_2 A \bigr\}.\nonumber
\end{eqnarray}
The event $\A_1^-$ is defined in the obvious related manner.
We may also define the one-arm event $\A_1$ on the ``uncolored'' space
$\HH$.
We will need the following extension of Lemma 2.4 in \cite{arXiv1008.1378} whose
proof applies easily to our present case.
The proof is analogous to that of Lemma 2.4 in \cite{arXiv1008.1378}, with the
difference that Theorem 5.8 and
Corollary 5.9 of \cite{\DuminilSixArm} replace the corresponding
results for Bernoulli percolation used in \cite{arXiv1008.1378}.

%
\begin{lemma}[(See Lemma 2.4 in \cite{arXiv1008.1378})]\label{l.lemma24}
Let $A$ be a piecewise smooth annulus in $[0,1]^2$. Then
\[
\P\bigl[\bar\omega_a \in\A_1^\pm\bigr] \longrightarrow
\P\bigl[\bar
\omega_\infty\in\A_1^\pm\bigr]
\]
as the mesh size $a\to0$.
Furthermore, in any coupling of the measures $\P_a$ and $\P_\infty$ on
the space $(\HH\times\HH, \T\otimes\T)$, in which $\bar\omega_a \to
\bar\omega_\infty$ a.s. we have that $\1_{\A_1^\pm}(\bar\omega_a) \to
\1_{\A_1^\pm}(\bar\omega_\infty)$ almost surely.
\end{lemma}

If $A$ is the annulus $[-c_1,c_1]^2 \setminus[-c_2/2,c_2/2]^2$, we define
$\alpha_1^\mathrm{FK}(c_2,c_1) = \Pb{\omega_\infty^+ \in\A_1}$.

\subsubsection{General setup of convergence: 
The space 
\texorpdfstring{$\mathscr{H}\times\mathscr{H}\times\mathcal{H}^{-3}$}
{HtimesHtimesH-3}}\label{ss.generalsetup}

Let us consider the coupling $(\bar\omega_a, \Phi^a(\bar\omega_a)) =
(\omega_a^+, \omega_a^-, \Phi^a) \in\HH\times\HH\times\Sob^{-3}$.
In order to prove our main Theorem~\ref{th.main}, will prove the
following stronger result.

%
\begin{theorem}[(Under Assumption~\ref{a.uniqueness})]\label{th.uniqField}
The random variables $(\bar\omega_{a}, \Phi^{a})\in\HH\times\HH
\times\Sob^{-3}$ converge in law as the mesh size $a\to0$
to $(\bar\omega_\infty, \Phi^\infty)$ for the topology induced by the metric\footnote{Recall from Section~\ref{ss.notations} that we have chosen a
metric $d_\HH$ on $\HH$ which induces the topology~$\T$.}
$d_\HH\oplus d_\HH\oplus\|\cdot\|_{\Sob^{-3}}$.

Furthermore, the limiting random variable $\Phi^\infty\in\Sob^{-3}$ is
measurable with respect to $\bar\omega_\infty$, that is, we have %
\[
\Phi^\infty= \Phi^\infty(\bar\omega_\infty).
\]
\end{theorem}

From Proposition~\ref{pr.tight}, we already know that $(\bar\omega_a,
\Phi^a(\bar\omega))$ is tight in the space $\HH\times\HH\times\Sob
^{-3}$ endowed with the metric $d_\HH\oplus d_\HH\oplus\| \cdot\|
_{\Sob^{-3}}$. As in Corollary~\ref{cor.sub}, one thus has
subsequential scaling limits: that is, one can find a subsequence $a_k
\to0$ such that $(\bar\omega_{a_k}, \Phi^{a_k})$ converges in law to
$(\bar\omega_\infty, \Phi^*)$ (here we use Assumption~\ref{a.uniqueness} which says that there is a unique possible subsequential
scaling limit for $\bar\omega_a$).
Since the space $(\HH\times\HH\times\Sob^{-3}, d_\HH\oplus d_\HH
\oplus\| \cdot\|_{\Sob^{-3}})$ is a complete separable metric space,
one can apply Skorohod's theorem. This gives us a joint coupling of the
above processes such that
%
%
\begin{equation}
\label{e.conv} \bigl(\bar\omega_{a_k}, \Phi^{a_k} \bigr)
\overset{\mathrm{a.s.}} { \longrightarrow} \bigl(\bar\omega_\infty, \Phi^* \bigr).
\end{equation}

Proving Theorem~\ref{th.uniqField} boils down to proving that $\Phi^*$
is in fact measurable with respect to $\bar\omega_\infty$.
Achieving this would indeed conclude the proof of Theorem~\ref{th.uniqField} (and thus Theorem~\ref{th.main}) since it would
uniquely characterize the subsequential scaling limits of $(\bar\omega
_a, \Phi^a)$.

The purpose of the next subsection is to reduce the proof of Theorem
\ref{th.uniqField} to the study of the renormalized
magnetization in a square box.

\subsubsection{Reduction to the renormalized magnetization 
in a dyadic box}

We wish to prove Theorem~\ref{th.uniqField}, that is, to show that if
$(\bar\omega_{a_k}, \Phi^{a_k}) \overset{\mathrm{a.s.}}{\longrightarrow} (\bar
\omega_\infty, \Phi^*)$,
then~$\Phi^*$ can be expressed as a measurable function of $\bar\omega
_\infty$. Since $\Phi^*\in\Sob^{-3}$,
it can be decomposed in the orthonormal basis $\{e_{j,k}\}_{j,k\geq
1}$, introduced\vspace*{1pt} in Appendix~\ref{tight}, of the
space $C^{\infty}_0([0,1]^2)$ endowed with the $L^2$ norm:
if $f = \sum_{j,k \geq1} a_{j,k} e_{j,k} \in C^{\infty}_0([0,1]^2)$, then
\[
\bigl\langle\Phi^*,f\bigr\rangle = \sum_{j,k\geq1}
\bigl\langle \Phi^*,e_{j,k}\bigr\rangle a_{j,k},
\]
which we write as
%
\begin{equation}
\label{e.decomp} \Phi^* = \sum_{j,k\geq1} \bigl\langle
\Phi^*,e_{j,k}\bigr\rangle e_{j,k}.
\end{equation}
From the a.s. convergence~\eqref{e.conv}, we have that for any fixed
$j,k\geq1$:
\[
  \bigl(\bar\omega_{a_k}, \bigl\langle\Phi^{a_k},
e_{j,k}\bigr\rangle \bigr) \overset{\mathrm{a.s.}} { \longrightarrow} \bigl(\bar
\omega_\infty, \bigl\langle\Phi^*, e_{j,k}\bigr\rangle \bigr),
\]
where the convergence holds for the metric $d_\HH\oplus d_\HH\oplus\|
\cdot\|_{\R}$.
In order to prove Theorem~\ref{th.uniqField}, thanks to the
decomposition \eqref{e.decomp}, we only need
to prove that for each fixed $j,k\geq1$, the limiting quantity $ \<
{\Phi^*, e_{j,k}}$ is itself measurable w.r.t. $\bar\omega_\infty$.

It turns out that one can further reduce the difficulty of this task by
approximating the functions $\{e_{j,k}\}_{j,k \geq1}$
using step functions as follows. Let us fix some $j,k\geq1$. For any
small $\beta>0$, one can find dyadic squares
$B_i$ and real numbers $b_i$ so that if $g_\beta:= \sum_i b_i 1_{B_i}$, then
\[
\| e_{j,k} - g_\beta\|_{L^\infty([0,1]^2)} < \beta.
\]

Now, exactly as in the proof of Lemma~\ref{l.SOB}, it is not hard to
check that
\[
  \E\bigl[ \bigl( \bigl\langle\Phi^a, e_{j,k}\bigr
\rangle - \bigl\langle \Phi^a, g_\beta\bigr\rangle
\bigr)^2\bigr] \le C \beta^2,
\]
uniformly in $a>0$ (for some universal constant $C>0$).
If one can show that, as $a_k \searrow0$, $\<{\Phi^{a_k}, g_\beta}$
converges a.s. to a measurable function $G_\beta(\bar\omega_\infty)$,
then using the uniform $L^2$ bounds from Appendix~\ref{tight} together
with the triangle inequality in~$L^2$, it follows that $G_\beta(\bar
\omega_\infty)$ converges as $\beta\to0$ in $L^2$ to $\<{\Phi
^*,e_{j,k}}$. Since $L^2$ is complete, $G_\beta(\bar\omega_\infty)$
has an $L^2$-limit $G_0$ as $\beta\to0$ which is itself measurable
w.r.t. $\bar\omega_\infty$ and one has necessarily that $\<{\Phi^*,
e_{j,k}}\overset{\mathrm{a.s.}}{=} G_0(\bar\omega_\infty)$.
Since $G_\beta$ is a linear combination of magnetizations in dyadic
squares, it follows from the above discussion that Theorem~\ref{th.uniqField} is a corollary of the following theorem.

%
\begin{theorem}\label{th.uniqMagn}
Let $B$ be any dyadic square in $[0,1]^2$ and let the renormalized
magnetization in $B$ be the random variable
\[
  m^a_B=m^a:= a^{15/8}
\sum_{x\in a \Z^2 \cap B} \sigma_x.
\]

Then the coupled random variable $(\bar\omega_{a}, m^{a})\in\HH\times
\HH\times L^2$ converges in law as the mesh size $a\to0$ to $(\bar
\omega_\infty, m)$
for the topology induced by the metric $d_\HH\times d_\HH\times\|
\cdot\|_{L^2}$.
Furthermore, the limiting random variable $m\in L^2$ is measurable with
respect to $\bar\omega_\infty$, that is, we have %
\[
m = m(\bar\omega_\infty).
\]
\end{theorem}

We now turn to the proof of this theorem. Without loss of generality
and for the sake of simplicity, we may assume that our dyadic square
$B$ is just $[0,1]^2$.

\subsection{Scaling limit for the magnetization random 
variable (proof of Theorem~\texorpdfstring{\protect\ref{th.uniqMagn}}{2.6})}

\subsubsection{Structure of the proof of Theorem \texorpdfstring{\protect\ref{th.uniqMagn}}{2.6}}

The setup for the scaling limit of~$m^a$ is similar to the setup we
explained above (in Section~\ref{ss.generalsetup})
for the scaling limit of $\Phi^a$. Namely, we consider the coupling
$(\bar\omega_a, m^a)$ embedded in the metric space $( \HH^2 \times L^2,
d_\HH\oplus d_\HH\oplus\|\cdot\|_{L^2}) $.
The tightness of $(\bar\omega_a, m^a)$ easily follows from the
stronger tightness of Proposition~\ref{pr.tight} (see also \cite
{\CamiaNewmanIsing} and \cite{Camia}). In particular, there exist
subsequential scaling limits
\[
\bigl(\bar\omega_{a_k}, m^{a_k} \bigr) \overset{d} {
\longrightarrow} \bigl(\bar\omega_\infty^*, m^* \bigr).
\]
By Assumption~\ref{a.uniqueness}, there is a unique possible law for
$\bar\omega_\infty^*$, which we denoted by~$\bar\omega_\infty$. In
order to prove Theorem~\ref{th.uniqMagn}, it remains to show that the
second coordinate $m^*\in L^2$ is measurable with respect to the first one.

As previously, let us couple all these random variables using
Skorohod's theorem so that
%
%
\begin{equation}
\label{e.Skoro} \bigl(\bar\omega_{a_k}, m^{a_k} \bigr)
\overset{\mathrm{a.s.}} {\longrightarrow} \bigl(\bar\omega_\infty, m^* \bigr)
\end{equation}
for the metric $d_\HH\oplus d_\HH\oplus\|\cdot\|_{L^2}$.

The main idea will be to approximate the quantity $m^a$ by relying only
on ``macroscopic information'' from the coupled configuration $\bar
\omega_a$. The ``macroscopic quantities'' we are allowed to use are the
quantities which are preserved in the scaling limit $\bar\omega_a \to
\bar\omega_\infty$ (i.e., crossing events and so on, see Section~\ref{ss.measurablevents}).

We will approximate the magnetization $m^a$ by a two step procedure.
Roughly speaking, we first approximate the magnetization as a rescaled
sum of spin variables $\sigma_x$ such that
$x$ is the starting point of a ``macroscopic'' FK path, and then
approximate the latter sum by means of
``macroscopic'' quantities following \cite{arXiv1008.1378}, as explained below.

For the first step, we fix some small dyadic scale $\rho\in\{2^{-k},
k\in\N\}$ and divide the square
$[0,1]^2$ along the grid $\rho\Z^2$. Let $S_\rho$ be the set of $\rho
$-squares thus obtained. For each
$\rho$-square $Q\in S_\rho$, consider the annulus $A_Q:=3Q\setminus Q$
where we denote by~$3Q$ the
square of side-length $3\rho$ centered on $Q$. We will divide the
clusters in the FK-configuration $\omega_a$
in two groups: the clusters which cross at least one annulus $A_Q, Q\in
S_\rho$ and the clusters which do not
cross any annulus. We may rewrite the magnetization $m^a$ as follows:
%
%
\begin{eqnarray}
  m^a &=& \sum_{x\in[0,1]^2 \cap a\Z^2}
a^{15/8} \sigma_x
\\
&=& \sum_{Q\in S_\rho} \biggl( \sum
_{x\in Q\dvtx x \leftrightarrow\partial(3Q)} a^{15/8} \sigma_x + \sum
_{x \in Q\dvtx x \nleftrightarrow\partial(3Q)} a^{15/8} \sigma_x \biggr).
\end{eqnarray}

Following \cite{\CamiaNewmanIsing} (with a slightly different setup
here), let us show that the contribution of the second
inside sum is negligible in $L^2$. Indeed,
%
%
\begin{eqnarray}
&&  \biggl\| \sum_{Q\in S_\rho} \sum
_{x \in Q\dvtx x \nleftrightarrow\partial(3Q)} a^{15/8} \sigma_x
\biggr\|_2^2\nonumber
\\
&&\qquad  =  \sum_{Q,Q'}\sum_{x,y} a^{15/4} \Eb{\sigma_x
\sigma_y 1_{x \in Q\dvtx x \nleftrightarrow\partial(3Q)} 1_{y \in
Q'\dvtx y
\nleftrightarrow\partial(3Q')}}
\nonumber
\\
&&\qquad = \sum_{Q,Q'}\sum_{x,y}
a^{15/4} \Eb{1_{x \leftrightarrow y} 1_{x \in Q\dvtx x \nleftrightarrow
\partial(3Q)} 1_{y \in Q'\dvtx y \nleftrightarrow\partial(3Q')}}
\\
&&\qquad \le \sum_{x,y\dvtx |x-y|\le8 \rho} a^{15/4}
\Eb{1_{x \leftrightarrow y}}
\\
&&\qquad  =  a^{15/4} O \biggl( a^{-2} \biggl(\frac{8 \rho} a
\biggr)^2 (\rho/a)^{-1/4} \biggr)
\\
&&\qquad = O \bigl(\rho^{7/4} \bigr).
\end{eqnarray}

Since we are looking for a limiting law for $m^a$ in $L^2$, it is thus
enough [up to a~small error of $O(\rho^{7/4})$]
to focus on the first summand
\[
m^a_\rho:= \sum_{Q\in S_\rho} \sum
_{x\in Q\dvtx x \leftrightarrow\partial(3Q)} a^{15/8} \sigma_x.
\]

Since $\rho>0$ is fixed and the mesh size $a\to0$, we are getting
closer to an approximation by ``macroscopic quantities.''
We still need to approximate in a suitable macroscopic manner the quantity
\[
L_Q^a= L_Q^a(\bar
\omega_a):=\sum_{x\in Q\dvtx x \leftrightarrow\partial(3Q)} a^{15/8}
\sigma_x
\]
for each $\rho$-square $Q\in S_\rho$. This is the second step of our
approximation procedure and for this
we will follow very closely the proof in \cite{arXiv1008.1378} of the scaling
limit of counting measures on pivotal points
(called \textit{pivotal measures}). In the rest of the proof, let us fix
the value of $\rho$ and fix some $\rho$-square $Q\in S_\rho$.

Let $\eps>0$ be some small fixed threshold (such that $a \ll\eps\ll
\rho$).
Divide the square $Q\cap a\Z^2$ into equal disjoint squares of
side-length $\eps_a:= a \floor{\eps/a}$. There are $N=\Omega(\eps
^{-2})$ such squares inside $Q$ (we do not need to keep the dependence
in $\rho$ in what follows) plus $O(\eps^{-1})$ squares which intersect
the boundary of $Q$. Let $ ( B_i )_{i\in\{1,\ldots, N\}}$ denote the
set of such $\eps_a$-squares inside $Q$.

For each $i\in[N]:= \{ 1, \ldots, N \}$, let
%
%
\begin{equation}
  X_i^\eps= X_i:= \sum
_{x\in B_i\dvtx x \leftrightarrow\partial(3Q)} a^{15/8} \sigma_x.
\end{equation}
Furthermore, let $B:=\cup B_i \subset Q$. We thus have
%
%
\begin{equation}
  L_Q^a = \sum
_{i\in[N]} X_i^\eps+ \sum
_{x\in Q\setminus B\dvtx x \leftrightarrow\partial(3Q)} a^{15/8} \sigma_x.
\end{equation}

The second term (which arises when $\rho$ is not a multiple of $\eps
_a$) turns out to be negligible in $L^2$ as well. Indeed,
%
%
\begin{eqnarray}
  \biggl\|\sum_{x\in Q\setminus B\dvtx x \leftrightarrow\partial(3Q)} a^{15/8}
\sigma_x\biggr\|_2^2 &\le& \sum
_{x,y \in Q\setminus B} a^{15/4} \Pb{x\leftrightarrow y}
\\
&\le& O(1) a^{15/4} \frac{\rho} {\eps} \sum_{k=1}^{\rho/\eps}
\biggl( \frac{\eps} a \biggr)^4 \biggl(\frac{a} {k \eps}
\biggr)^{1/4}
\\
& \le& O(1) \frac{\rho} {\eps} (\rho/\eps)^{3/4} \eps^{4}
\eps^{-1/4}
\\
&\le& O(1) \rho^{7/4} \eps^2.
\end{eqnarray}
%
Therefore, as $\eps$ goes to zero, and uniformly in the mesh size $a\le
\eps$,
the boundary term is negligible in $L^2$. It thus remains to control
the term $\sum_{i\in[N]} X_i^\eps$.

For this, let us introduce for each $i\in[N]$, the variables
%
%
\begin{equation}
\label{y-function} Y_i^\eps=Y_i^\eps(
\bar\omega_a):= 1_{\{ B_i \overset{\omega_a^+}{\longleftrightarrow}
\partial(3Q)\} } - 1_{\{B_i \overset{\omega_a^-} {\longleftrightarrow}
\partial(3Q)\}}.
\end{equation}
We will prove in the next subsection the following proposition.

\begin{proposition}\label{pr.XY}
There exists a universal constant $c>0$ such that for any square $Q$ of
side-length $\rho$ as above,
we have
%
%
\begin{equation}
\label{e.XY} \biggl\| \sum_{i \in[N]} X_i^\eps
- c \beta(\eps) \sum_{i \in[N]} Y_i^\eps
\biggr\|_2 \longrightarrow0
\end{equation}
as $\eps\to0$ uniformly in $a\le\eps$ and where $\beta(\eps):= \eps^2
\alpha_1^\mathrm{FK}(\eps,1)^{-1}$
and $\alpha_1^\mathrm{FK}(\eps,1)$ is the probability of the one-arm
event in the annulus
$[-1,1]^2 \setminus[-\eps/2,\eps/2]^2$ for $\omega_\infty^+$.
\end{proposition}

Before proving the proposition, let us explain why it indeed implies
Theorem~\ref{th.uniqMagn}.
From Section~\ref{ss.measurablevents},
it follows that the functions $Y_i^\eps$ defined in \eqref{y-function}
can be seen as measurable functions
of $\bar\omega_\infty$ and that, for each $i\in[N]$, along the above
subsequence $(a_k)$, one has
[see equation (\ref{e.Skoro})]:
\[
  Y_i^\eps(\bar\omega_{a_k})
\overset{\mathrm{a.s.}} {\longrightarrow} Y_i^\eps(\bar
\omega_\infty).
\]

Furthermore, one can see from Proposition~\ref{pr.XY} that $\| \beta
(\eps) \sum Y_i^\eps(\bar\omega_a) \|_2$
is bounded uniformly in $0<a\le\eps$. This implies, modulo some
triangle inequalities, that
\[
  \Bigl\| L_Q^a(\bar\omega_{a_k}) - c
\beta(\eps) \sum Y_i^\eps(\bar
\omega_\infty) \Bigr\|_2 \longrightarrow0,
\]
uniformly in $0<a_k<\eps$. This in turn implies that the sequence $ ( c
\beta(\eps)\times\break  \sum Y_i^\eps(\bar\omega_\infty) )_{\eps>0}$ is a
Cauchy-sequence in $L^2$. In particular, it has an $L^2$-limit that we
may denote by $L_Q(\bar\omega_\infty)$ and this $L^2$-limit is such that
\[
  \bigl\| L_Q^{a_k}(\bar\omega_{a_k}) -
L_Q(\bar\omega_\infty) \bigr\|_2 \longrightarrow0
\]
as the mesh size $a_k\to0$.

Using the above estimates, we have that
\[
  \biggl\|m^{a_k} - \sum_Q
L_Q(\bar\omega_\infty) \biggr\|_2^2
\longrightarrow0,
\]
uniformly in $0<a_k<\rho$.
Exactly as above with the second order approximation in $\eps$, the
above displayed equation (plus the $L^2$ bounds we already have)
implies that the Cauchy sequence $ ( \sum_Q L_Q(\bar\omega_\infty)
)_{\rho>0}$ has an $L^2$-limit denoted by $m(\bar\omega_\infty)$ as
$\rho\to0$. Finally, thanks to the a.s. convergence in equation
(\ref{e.Skoro}), this $L^2$-limit must be such that
%
%
\begin{equation}
  m^* \overset{\mathrm{a.s.}} {=} m(\bar\omega_\infty),
\end{equation}
which completes the proof of Theorem~\ref{th.uniqMagn}, modulo proving
Proposition~\ref{pr.XY}.

\subsubsection{Proof of Proposition \texorpdfstring{\protect\ref{pr.XY}}{2.7}}

We want to show that for any $\delta>0$, one can take $\eps>0$
sufficiently small so that for any
$0<a<\eps< \rho$,
\[
\E\Bigl[ \Bigl(\sum X_i^\eps- c \beta(\eps) \sum
Y_i^\eps \Bigr)^2 \Bigr] \le\delta.
\]

Let us decompose this quantity as follows.
%
%
\begin{eqnarray}\label{e.jacques}
&& \E \Bigl[ \Bigl(\sum X_i^\eps- c
\beta(\eps) \sum Y_i^\eps \Bigr)^2 \Bigr]\nonumber
\\
&&\qquad  = \sum_{i,j} \E\bigl[ \bigl( X_i^\eps
- c \beta(\eps) Y_i^\eps \bigr) \bigl(
X_j^\eps- c \beta(\eps) Y_j^\eps
\bigr)\bigr]
\nonumber\\[-8pt]\\[-8pt]
&&\qquad \le\sum_{i,j\dvtx d(B_i,B_j)\le r} \bigl( \E\bigl[X_i^\eps
X_j^\eps\bigr] + c^2 \beta(\eps)^2
\E\bigl[Y_i^\eps Y_j^\eps\bigr] \bigr)\nonumber
\\
&&\quad\qquad {}+ \sum_{i,j\dvtx d(B_i, B_j)> r} \E \bigl[ \bigl( X_i^\eps
- c \beta(\eps) Y_i^\eps \bigr) \bigl(
X_j^\eps- c \beta(\eps) Y_j^\eps
\bigr)\bigr],\nonumber
\end{eqnarray}
where $r$ is a mesoscopic scale $\eps\ll r \ll\rho$ which will be
chosen later. To go from the first to the second line, we used that
fact that the cross product terms are necessarily negative as can be
seen by first conditioning on the noncolored FK configuration $\omega_a$.

The first term of the RHS of the above displayed inequality is easy to
bound. Indeed,
\begin{eqnarray*}
\sum_{i,j\dvtx d(B_i,B_j)\le r} \E\bigl[X_i^\eps
X_j^\eps\bigr] & \le&\sum_{x,y\in a \Z^2 \cap Q\ \mathrm{s.t.}\ d(x,y)\le2 r}
a^{15/4} \Eb{\sigma_x \sigma_y} \le O \bigl(
r^{7/4} \bigr)
\end{eqnarray*}
and similarly for $\sum_{i,j\dvtx d(B_i,B_j)\le r} c^2 \beta(\eps)^2
\Eb
{Y_i^\eps Y_j^\eps}$.
One can thus fix $r>0$ small enough so that, uniformly in $a<\eps<r$,
the first term in the RHS
of \eqref{e.jacques} is~$<\delta/2$.

For the second term, we proceed as in \cite{arXiv1008.1378} using a coupling argument.
Proposition~\ref{pr.XY} will follow from the next lemma.

\begin{lemma}\label{l.ij}
For any fixed $r<\rho< 1$ and any $\tilde\delta>0$, one can choose
$\eps=\eps(r,\rho, \tilde\delta)>0$ small enough such that for any
pair of squares $B_i, B_j$ with $d(B_i,B_j)>r$, one has
\[
  \E \bigl[ \bigl( X_i^\eps- c \beta(\eps)
Y_i^\eps \bigr) \bigl( X_j^\eps- c
\beta(\eps) Y_j^\eps \bigr)\bigr] \le\frac{\tilde\delta} 2
\E \bigl[X_i^\eps X_j^\eps\bigr].
\]
\end{lemma}

Let us explain why this lemma is enough to conclude the proof.
Summing the estimate provided by the lemma over all $B_i, B_j$ with
$d(B_i,B_j)>r$, one gets
\begin{eqnarray*}
  \sum_{i,j\dvtx d(B_i, B_j)> r} \E \bigl[ \bigl(
X_i^\eps- c \beta(\eps) Y_i^\eps
\bigr) \bigl( X_j^\eps- c \beta(\eps) Y_j^\eps
\bigr)\bigr] &\le&\frac{\tilde\delta} 2 \E \bigl[ \bigl(L_Q^a(\bar
\omega_a) \bigr)^2\bigr].
\end{eqnarray*}
Now, it is straightforward to check that the second moment $\Eb
{(L_Q^a(\bar\omega_a))^2}$ is bounded by $C \rho^{15/4}$ uniformly in
$a<\eps<\rho$ where $C$ is some universal constant. By choosing $\tilde
\delta= \delta/C$, we conclude the proof.

\begin{pf*}{Proof of Lemma~\ref{l.ij}}
Let us fix two squares $B_i$ and $B_j$ at distance at least~$r$ from
each other.
Conditioned on the event that both $B_i$ and $B_j$ are connected to~$\partial(3Q)$, our strategy
is to compare how things look within the $\eps$-square $B_i$ with the
following ``test case.''
Consider the $\eps$-square $B_0$ centered at the origin and let $Q_0$
be the square
$[-\rho,\rho]^2$ also centered at the origin. Let us define
%
%
\begin{equation}
  X_0:= a^{15/8} \sum
_
{
x\in a\Z^2 \cap B_0
} \sigma_x 1_{\{ x \leftrightarrow\partial Q_0\ \mathrm{in}\ \bar
\omega_a \} }.
\end{equation}
Recall the events $\A_1^\pm$ defined in Section~\ref{ss.measurablevents} and applied here to the annulus
$A=Q_0 \setminus B_0$. We first wish to show that there is a constant
$c>0$ such that
%
%
\begin{equation}
  \cases{ \E \bigl[X_0 \md{\mathcal A}_1^+\bigr] \sim
c \beta(\eps),
\vspace*{3pt}\cr
\E\bigl[X_0 \md{\mathcal A}_1^-\bigr] \sim-c
\beta(\eps)}
\end{equation}
uniformly as $0<a<\eps$ go to 0.
To see why this holds we note that, as in Section~4.5 in \cite{arXiv1008.1378},
one has that
\[
  \E \bigl[X_0 \md{\mathcal A}_1^+\bigr] \sim(
\eps/a)^2 a^{15/8} \frac{\alpha_1^\mathrm{FK}(a, 1)} {
\alpha_1^\mathrm{FK}(\eps, 1) }.
\]
To adapt the proof from \cite{arXiv1008.1378}, it is enough to have
bounds on the half-plane exponents of critical FK percolation ``in the bulk.''
Such bounds follow from standard percolation arguments, using the RSW
theorem of \cite{arXiv0912.4253}.

Now, using Theorem 1.3 (with $k=0$) in \cite{arXiv1202.2838} together with
Wu's result, Theorem~\ref{th.Wu},
we have that $\alpha_1^\mathrm{FK}(a,1)\sim c a^{1/8}$ as $a\to0$,
which explains the desired asymptotic.

In what follows, for any $u\geq\eps$, we will denote by $D_u$
($\widetilde D_u$) the square centered around $B_i$ ($B_j$) of
side-length $u$.
Let\vspace*{1pt} us fix yet another mesoscopic scale $\gamma$ so that $\eps\ll
\gamma\ll r$ (e.g., $\gamma:=r^2$).
Let $m:=d(B_i,B_j)/2$ and let\vspace*{-2pt} $z$ be the midpoint between the centers
of $B_i$ and $B_j$.
Let $R^+$ be the event that $\{ \partial D_\gamma\overset{\omega
^+}{\longleftrightarrow} \partial D_{m} \} \cap\{\mbox{there is a
circuit of } \omega^+ \mbox{ inside } D_r \setminus D_\gamma
\mbox{ that surrounds } D_\gamma\}$ (see Figure~\ref{fig:R+}). The event $R^-$~is defined
similarly. Notice that $R^+ \cap R^- = \varnothing$. On the event $R^\pm
$, let $C=C(\bar\omega)$ be the outermost such open circuit for the FK
configuration
$\omega\in\HH$ (the outermost open circuit necessarily has the
appropriate color).

%
\begin{figure}

\includegraphics{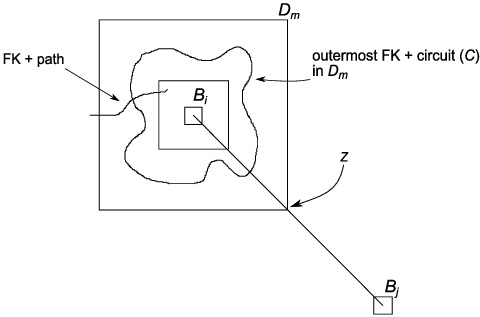}

\caption{The event $R^{+}$.}\label{fig:R+}
\end{figure}

Let us analyze in the term $ \Eb{ ( X_i - c \beta(\eps) Y_i^\eps) (
X_j - c \beta(\eps) Y_j ^\eps)}$
the contribution coming from the event $(R^+ \cup R^-)^c$, namely,
\begin{eqnarray*}
&& \E \bigl[ \bigl( X_i - c \beta(\eps) Y_i^\eps
\bigr) \bigl( X_j - c \beta(\eps) Y_j ^\eps
\bigr); \bigl(R^+ \cup R^- \bigr)^c\bigr]
\\
&&\qquad \le \E \bigl[X_i
X_j + c^2 \beta(\eps)^2 Y_i
Y_j; \bigl(R^+ \cup R^- \bigr)^c\bigr].
\end{eqnarray*}
See the explanation after~\eqref{e.jacques} as to why the cross product
terms are negative.
Following \cite{arXiv1008.1378},
\begin{eqnarray*}
&& \E \bigl[X_iX_j; \bigl(R^+ \cup R^-
\bigr)^c\bigr]
\\
&&\qquad= \E \bigl[X_i X_j; \partial D_\eps
\overset{ \omega} {\longleftrightarrow} \partial D_m;
\\
&&\hspace*{44pt} \partial
\widetilde D_\eps\overset{\omega} {\longleftrightarrow} \partial
\widetilde D_m;\partial B(z, 2m) \overset{\omega} {
\longleftrightarrow} \partial(3Q); \bigl(R^+ \cup R^- \bigr)^c\bigr]
\\
&&\qquad \le O(1) \EFK{} {\mathrm{wired}} {\widetilde X_i} \EFK{} {
\mathrm{wired}} {\widetilde X_j} \alpha_1^{\mathrm{wired}}(
\eps, \rho) \alpha_1^{\mathrm{wired}}(\eps, \gamma)
\\
&&\quad\qquad{}\times  \P \bigl[\partial
D_\gamma\overset{\omega} {\longleftrightarrow} \partial
D_m; \bigl(R^+ \cup R^- \bigr)^c\bigr],
\end{eqnarray*}
where we have just used FKG and where we dominated $X_i$ by $\widetilde
X_i$, the number of points in $B_i$ connected to $\partial B_i$ (we
also used some straightforward quasi-multiplicativity for the one-arm
FK event which follows easily from the RSW theorem in \cite{arXiv0912.4253};
see, e.g., \cite{\Wer09} for an explanation of
quasi-multiplicativity in the case of standard percolation and see \cite
{\DuminilSixArm} for quasi-multiplicativity results in the case of FK
percolation). Now, using FKG with RSW from \cite{arXiv0912.4253}, we get that
there exists an exponent $\xi>0$ such that
$\Pb{R^+ \cup R^- \md\partial D_\gamma\overset{\omega
}{\longleftrightarrow} \partial D_m} \geq1- (\gamma/ m)^\xi$, which
implies that
\[
  \P \bigl[\partial D_\gamma\overset{\omega} {
\longleftrightarrow} \partial D_m; \bigl(R^+ \cup R^-
\bigr)^c\bigr] \le\alpha_1^\mathrm{FK}(\gamma, m) (
\gamma/m)^\xi.
\]

Altogether, we obtain that
\[
  \E \bigl[X_iX_j; \bigl(R^+ \cup R^-
\bigr)^c\bigr] \le O(1) \Eb{X_i X_j} (\gamma/
m)^\xi.
\]
The term $\Eb{c^2 \beta(\eps)^2 Y_i Y_j; (R^+ \cup R^-)^c}$ can be
treated similarly.
We may thus focus our analysis on what is happening on the event $R^+
\cup R^-$.
Let $\mathcal{F}_C$ be the filtration induced by the configuration
outside the contour $C$.
One can write
\begin{eqnarray*}
&&\bigl|\E \bigl[ \bigl( X_i - c \beta(\eps) Y_i^\eps
\bigr) \bigl( X_j - c \beta(\eps) Y_j ^\eps
\bigr); R^+ \cup R^-\bigr]\bigr|
\\
&&\qquad = \E \bigl[\bigl|X_j - c\beta(\eps) Y_j\bigr|
\E \bigl[\bigl|X_i - c \beta(\eps) Y_i\bigr| \md\mathcal{F}_C\bigr];
R^+ \cup R^-\bigr],
\end{eqnarray*}
since on the event $R^+ \cup R^-$, the variable $|X_j - c\beta(\eps)
Y_j|$ is measurable w.r.t. $\mathcal{F}_C$.
Now,
\begin{eqnarray*}
&& \E \bigl[\bigl|X_j - c\beta(\eps) Y_j\bigr| \E \bigl[\bigl|X_i
- c \beta(\eps) Y_i\bigr| \md\mathcal{F}_C\bigr]; R^+ \cup R^-\bigr]
\\
&&\qquad = \P \bigl[R^+\bigr] \E \bigl[\bigl|X_j - c\beta(\eps) Y_j\bigr|
\E \bigl[\bigl|X_i - c \beta(\eps) Y_i\bigr| \md\mathcal{F}_C\bigr]
\md R^+ \bigr]
\\
&&\quad\qquad{} + \P \bigl[R^-\bigr] \E \bigl[\bigl|X_j - c\beta(\eps)
Y_j\bigr| \E \bigl[\bigl|X_i - c \beta(\eps) Y_i\bigr| \md
\mathcal{F}_C\bigr] \md R^- \bigr].
\end{eqnarray*}

Let us analyze the first term, it gives
\begin{eqnarray*}
&& \P \bigl[R^+\bigr] \E \bigl[\bigl|X_j - c\beta(\eps) Y_j\bigr|
\E \bigl[\bigl|X_i - c \beta(\eps) Y_i\bigr| \md\mathcal{F}_C,
R^+\bigr] \md R^+ \bigr]
\\
&&\qquad = \P \bigl[R^+\bigr] \E \bigl[\bigl|X_j - c\beta(\eps) Y_j\bigr|
\1_{C \leftrightarrow\partial(3Q)} \E \bigl[\bigl|X_i - c \beta(\eps) Y_i\bigr| \md C,
R^+\bigr] \md R^+ \bigr]
\\
&&\qquad = \P \bigl[R^+\bigr] \E \bigl[\bigl|X_j - c\beta(\eps) Y_j\bigr|
\1_{C \leftrightarrow\partial(3Q)} \P \bigl[\partial B_i \leftrightarrow C \md C, R^+\bigr]
\\
&&\hspace*{103pt}{}\times \E \bigl[\bigl|X_i - c \beta(\eps) Y_i\bigr| \md C, \partial
B_i \leftrightarrow C\bigr] \md R^+ \bigr].
\end{eqnarray*}

We will prove below the following lemma.

\begin{lemma}[(Coupling lemma)]\label{l.arg}
For any contour $C$, we have the following control on the conditional
expectation:
\[
  \E \bigl[\bigl|X_i - c \beta(\eps) Y_i\bigr| \md C,
\partial B_i \leftrightarrow C\bigr] \le K (\eps/\gamma)^\alpha
\beta(\eps) 
\]
for some exponent $\alpha>0$ and some constant $K \in(0,\infty)$.
\end{lemma}

Plugging this lemma into the last displayed equation leads to
%
%
\begin{eqnarray}
\label{e.cusp} && \P \bigl[R^+\bigr] \E \bigl[\bigl|X_j - c\beta(\eps)
Y_j\bigr| \E \bigl[\bigl|X_i - c \beta(\eps) Y_i\bigr| \md
\mathcal{F}_C, R^+\bigr] \md R^+ \bigr]
\nonumber
\\
&&\qquad \le C(\eps/\gamma)^\alpha\beta(\eps) \P \bigl[R^+\bigr]\nonumber
\\
&&\quad\qquad {}\times
\E \bigl[\bigl|X_j - c\beta(\eps) Y_j\bigr| \1_{C \leftrightarrow\partial3Q} \P \bigl[
\partial B_i \leftrightarrow C \md C, R^+\bigr] \md R^+ \bigr]
\\
&&\qquad \le O(1) (\eps/\gamma)^\alpha\beta(\eps) \alpha_1^\mathrm{FK}(
\gamma, m)\nonumber
\\
&&\quad\qquad{}\times  \E \bigl[\bigl|X_j - c\beta( \eps) Y_j\bigr|
\1_{C \leftrightarrow\partial(3Q)} \P \bigl[\partial B_i \leftrightarrow C \md C, R^+\bigr] \md
R^+ \bigr].\nonumber
\end{eqnarray}

Now, similarly to the above analysis of what happens on the event $(R^+
\cup R^-)^c$, it is not hard to check
by cutting into different scales and dominating by wired boundary
conditions that
\begin{eqnarray*}
&&  \E \bigl[\bigl|X_j - c\beta(\eps) Y_j\bigr|
\1_{C \leftrightarrow\partial(3Q)} \P\bigl[\partial B_i \leftrightarrow C \md C, R^+\bigr] \md
R^+ \bigr]
\\
&&\qquad  \le O(1) \beta(\eps) \alpha_1^\mathrm{FK}( \eps, \rho)
\alpha_1^\mathrm{FK}(\eps, \gamma),
\end{eqnarray*}
which together with~\eqref{e.cusp} and quasi-multiplicativity gives us
(since one has also the same estimate on the event $R^-$),
\begin{eqnarray*}
&&  \bigl| \E\bigl[ \bigl( X_i - c \beta(\eps) Y_i^\eps
\bigr) \bigl( X_j - c \beta(\eps) Y_j ^\eps
\bigr); R^+ \cup R^-\bigr]\bigr|
\\
&&\qquad \le  O(1) (\eps/\gamma)^\alpha\beta(
\eps)^2 \frac{\alpha_1^\mathrm{FK}(\eps, \rho)^2} {\alpha_1^\mathrm
{FK}(m,\rho)}
\\
&&\qquad \le O(1) (\eps/\gamma)^\alpha\Eb{X_i X_j},
\end{eqnarray*}
which (modulo proving Lemma~\ref{l.arg}) completes our proof of Lemma
\ref{l.ij}.
\end{pf*}


%


\subsubsection{Proof of Lemma \texorpdfstring{\protect\ref{l.arg}}{2.9}}

Let $\nu_C$ be the wired FK probability measure conditioned on
$\partial B_i\leftrightarrow C$ and let $\nu_0$ be the FK probability
measure in $Q_0$ conditioned on the event $\A_1=\A_1(\tilde x + Q_0
\setminus B_0)$, where we translated the annulus $A = Q_0 \setminus
B_0$ so that it surrounds $B_i$. Clearly, in the domain $\mathcal{D}_C$
(inside the circuit $C$), the measure $\nu_C$ dominates $\nu_0$. Using
RSW from \cite{arXiv0912.4253}, there is an open circuit in $\mathcal{D}_C
\setminus B_i$ for $\omega^0_a\sim\nu_0$ with $\nu_0$-probability\vspace*{1pt} at
least $1- c(\eps/\gamma)^\xi$. Let us call this event $W$.
On the event $W$, let $\widetilde C$ be the outermost circuit inside
$\mathcal{D}_C$ for $\omega_a^0$. Since $\nu_C$ dominates $\nu_0$, one
can couple $\omega^C_a\sim\nu_C$ with $\omega^0_a$ so that on the
event $W$, they share the same open circuit $\widetilde C$ and are
conditioned inside $\mathcal{D}_{\widetilde C}$ only on
the constraint $\{ \partial B_i \leftrightarrow\widetilde C \}$; in
particular, on the event $W$, in this coupling one has $X_i = X_0$. In
order to prove Lemma~\ref{l.arg}, it is enough to show that $\Eb{X_i;
W^c \md C, \partial B_i \leftrightarrow C}$ and $\Eb{X_0; W^c \md C,
\partial B_i \leftrightarrow C}$ are negligible w.r.t. $\beta(\eps)$,
which is straightforward using the quantity $\Eb{\widetilde X_i}\times
\Pb{W^c}$ as we did previously while analyzing what happened on the
event $(R^+\cup R^-)^c$.

\section{Second proof of the scaling limit of 
\texorpdfstring{$\Phi^a$}{Phia} using the $n$-point functions of 
Chelkak, Hongler and Izyurov}\label{second}

In this part, we will give a different proof of Theorem~\ref{th.main},
using the recent breakthrough results of Chelkak, Hongler and Izyurov
in \cite{arXiv1202.2838}. From our tightness result obtained in Appendix~\ref{tight}, recall that there exist subsequential scaling limits $\Phi
^\star= \lim\Phi^{a_n}$ for the convergence in law in the space $\Sob^{-3}$.
We wish to prove that there is a unique such subsequential scaling limit.
For this, we will use the following classical fact (see, e.g., \cite
{\LedouxTalagrand}).

%
\begin{proposition}
If $h$ is a random distribution in $\Sob^{-3}$ (for the sigma-field
generated by the topology of $\|\cdot\|_{\Sob^{-3}}$),
then the law of $h$ is uniquely characterized~by
\[
\phi_h(f):=\E \bigl[e^{i \<{h,f}}\bigr]
\]
as a function of $f\in\Sob^3$.
\end{proposition}


Using the tightness property proved in Appendix~\ref{tight}, Theorem~\ref{th.main} will thus follow from the next result.

%
\begin{proposition}\label{pr.chf}
For any $f\in\Sob^3$, the quantity
\[
\phi_{\Phi^a}(f) = \E \bigl[e^{i \<{\Phi^a, f}}\bigr]
\]
converges as the mesh size $a\searrow0$.
\end{proposition}

The proof of this proposition will be divided into two main steps as follows:
\begin{enumerate}[1.]
\item[1.] First, we will show that $\Phi^a$ has ``uniform exponential
moments'' which will allow us to express its characteristic function using
\[
\phi_{\Phi^a}(f) = \E \bigl[e^{i \<{\Phi^a, f}}\bigr] = 1+ \sum
_{k\geq1} \frac{i^k \Eb{\langle\Phi^a, f\rangle^k}}{k !}.
\]

\item[2.] Then it remains to compute each $k$th moment $\Eb{\<{\Phi^a,
f}^k}$, that is, to show uniqueness as $a \to0$. For this, one uses
the scaling limit results from \cite{arXiv1202.2838} together with Proposition
\ref{pr.close} below which takes care of $k$-tuples of points in the
plane where at least two points are close to each other.
\end{enumerate}

Let us now state the main result we will use from \cite{arXiv1202.2838}.

%
\begin{theorem}[(\cite{arXiv1202.2838}, Theorem 1.3)]\label{th.hongler}
Let $\Omega$ be a bounded simply connected domain, and $\Omega_a$ be
discretizations of $\Omega$ (built from $\Omega\cap a \Z^2$).
We denote by $\xi$ the boundary conditions chosen on $\Omega$, and we
assume $\xi$ to be either $+$ or $\mathrm{free}$ here.
Then, for any $k\geq1$, there exist \emph{$k$-point functions}
\[
z_1,\ldots,z_k \in\Omega^k \mapsto \langle
\sigma_{z_1}, \ldots,\sigma_{z_k}\rangle_{\Omega}^\xi,
\]
so that for any $\eps>0$, as the mesh size $a\to0$ and uniformly over
all $z_1,\ldots, z_k\in\Omega$ at distance at least $\eps$ from
$\partial\Omega$ and from each other, one has
%
%
\begin{equation}
  \varrho(a)^{-k/2} \cdot\EFK{\Omega_a} {\xi} {
\sigma_{z_1}, \ldots,\sigma_{z_k}} \longrightarrow \langle
\sigma_{z_1},\ldots,\sigma _{z_k}\rangle_{\Omega}^\xi.
\end{equation}
[Recall that $\varrho(a)$ is the renormalization factor defined
in~\eqref{e.varrho}.]

Furthermore, the functions $ \<{\sigma_{z_1}, \ldots,\sigma
_{z_k}}_{\Omega}^\xi$ are conformally covariant in the following sense:
if $\phi\dvtx \Omega\to\Omega'$ is a conformal map, then
\[
\langle\sigma_{z_1},\ldots,\sigma_{z_k}\rangle_{\Omega}^\xi
= \langle\sigma_{\phi(z_1)},\ldots,\sigma_{\phi(z_k)}
\rangle_{\Omega'}^\xi \prod\bigl|\phi'(z_i)\bigr|^{1/8}.
\]
\end{theorem}

%
\begin{remark}
\begin{itemize}
\item It is noted in \cite{arXiv1202.2838} that although their Theorem~1.3 is
stated only for plus boundary conditions, the conclusions are
valid for free and other boundary conditions as well.
\item In \cite{arXiv1202.2838}, the discretization is slightly different,
which means that our $k$-point function
$\<{\sigma_{z_1}, \ldots, \sigma_{z_k}}_{\Omega}^\xi$ is equal to the
one of \cite{arXiv1202.2838} only up to a constant factor.
%
\item In most of this paper, we assume Wu's result, Theorem~\ref{th.Wu}. In particular, one can then use the above theorem with
$a^{k/8}$ instead of $\varrho(a)^{k/2}$ (and with yet a further
change of the $k$ point function by another scalar). See Section~\ref{Theta} for the analysis when one does not wish to assume Wu's result.
\end{itemize}
\end{remark}

\subsection{Exponential moments for the magnetization random variable}
In this section, we shall show that if $m^a$ denotes the magnetization
random variable $\<{\Phi^a, 1_{[0,1]^2}}$ (for wired or free boundary
conditions on the square $[0,1]^2$), then $m^a$ has exponential
moments. More precisely, we will prove the following.

\begin{proposition}\label{pr.expomoments}
For any $t\in\R$, and for any boundary condition $\xi$ on~$[0,1]^2$,
one has
\[
\limsup_{a \searrow0} \E^{\xi}\bigl[e^{t m^a}\bigr] < \infty.
\]
\end{proposition}

There are a number of ways to prove this proposition.
We present one based on the 
Griffiths--Hurst--Sherman inequality
from \cite{\GHS}. Let us state it here.

%
\begin{theorem}[(GHS inequality, \cite{\GHS})]\label{th.GHS}
Let $G=(V,E)$ be a finite graph. Consider a pair ferromagnetic Ising
model on this graph (i.e., the interactions $J_{ij}$ between vertices
$i\sim j$ are nonnegative)
and assume furthermore that the external field $\mathbf{h} =
(h_v)_{v\in V}$ (which may vary from one vertex to another) is also
nonnegative. Under these general assumptions, one has for any vertices
$i,j,k\in V$,
\[
\langle\sigma_i\sigma_j \sigma_k\rangle-
\bigl( \langle\sigma_i\rangle \langle\sigma_j
\sigma_k\rangle+\langle\sigma_j\rangle \langle
\sigma_i \sigma_k\rangle+\langle\sigma_k
\rangle \langle\sigma_i \sigma_j\rangle \bigr) +2
\langle\sigma_i\rangle \langle\sigma_j\rangle \langle
\sigma_k\rangle\leq0.
\]
\end{theorem}

This inequality has the following useful corollary (see, e.g., \cite{CGNNC}).

\begin{corollary}\label{cor.GHS}
Let $G=(V,E)$ be a finite graph and let $K \subset V$ be a nonempty
subset of the vertices.
Let us consider a ferromagnetic Ising model on $G$ with the spins in
$K$ prescribed to be
$+$ spins and with a constant magnetic field $h\geq0$ on $V\setminus
K$. Then the partition
function of this model, that is,
\[
Z_{\beta, h}:= \sum_{\sigma\in\{-,+\}^{V\setminus K}} \exp{ \biggl(-\beta
E(\sigma) + h \sum_{i \in V \setminus K} \sigma_i
\biggr)},
\]
where $E(\sigma) = \sum_{i \sim j \in V} J_{ij} \sigma_i \sigma_j$, satisfies
\[
\partial^3_h \log(Z_{\beta,h}) \leq0.
\]
\end{corollary}

\begin{pf*}{Proof of Proposition~\ref{pr.expomoments}}
If $t<0$, using the symmetry of the Ising model, by changing the
boundary condition $\xi$ into $-\xi$,
we can assume $t>0$.
Hence, one may assume that $t\geq0$. This makes the function $x\mapsto
e^{t x}$ increasing, and one can thus use
the FKG inequality which implies that for any $t\geq0$ and any
boundary condition $\xi$, one has
\[
\E^{\xi} \bigl[e^{t m^a}\bigr] \le\E^{+} \bigl[e^{t m^a}\bigr].
\]
With $+$ boundary condition on $[0,1]^2\cap a\Z^2$, one can now rely on
the above corollary of the GHS inequality
which yields
\begin{eqnarray*}
\partial_h^3 \Bigl[ \log \Bigl(\sum
e^{\beta E(\sigma) + h \sum\sigma_i} \Bigr) \Bigr] & =& \partial_h^3
\biggl[ \log \biggl( \frac{ \sum e^{\beta E(\sigma) + h \sum\sigma_i} }{\sum e^{\beta E(\sigma)}} \biggr) \biggr]
\\
&=& \partial_h^3 \bigl[ \log\E_{\beta}\bigl[e^{h \sum\sigma_i}\bigr] \bigr]
\\
& \le& 0.
\end{eqnarray*}

With $\beta=\beta_c$ and $h:= t a^{15/8}$, one obtains
that for any $t\geq0$ and any mesh size $a>0$:
%
%
\begin{equation}
\label{e.thirdder} \partial_t^3 \log\E^{+}
\bigl[e^{t m^a}\bigr] \leq0.
\end{equation}

Now let $\phi(t):= \EFK{}{+}{e^{t m^a}}$. It is easy to check that
\[
\cases{ \phi'(0) = \E^{+}\bigl[m^a\bigr],
\vspace*{3pt}\cr
\phi''(0) = \E^{+}\bigl[ \bigl(m^a - \bigl
\langle m^a \bigr\rangle \bigr)^2\bigr].}
\]
This, together with~\eqref{e.thirdder}, implies that for any $t\geq0, a>0$:
\[
\log\E^{+}\bigl[e^{t m^a}\bigr] \le t \E^{+}\bigl[m^a\bigr] +
\frac{t^2} 2 \E^{+}\bigl[ \bigl(m^a - \bigl\langle
m^a \bigr\rangle \bigr)^2\bigr].
\]
By our choice of rescaling, $m^a:= a^{15/8} \sum_{x\in[0,1]^2\cap a \Z
^2}\sigma_x$, we know from
Proposition~\ref{pr.1and2} in Appendix~\ref{s.appendix} that
$\sup_{a>0} t \EFK{}{+}{m^a} + \frac{t^2} 2 \EFK{}{+}{(m^a - \langle
m^a\rangle)^2} = O(t+t^2) <\infty$,
which completes the proof of Proposition~\ref{pr.expomoments}.
\end{pf*}

We have the following easy corollary of Proposition~\ref{pr.expomoments}; it applies, for example, to $\Omega= [0,1]^2$ with
$\xi$ plus or free, where there is a unique limit, and also to quite
general $\Omega$ and $\xi$ where there may only
be limits along subsequences of $a \to0$.

%
\begin{corollary}
If $m=m^\xi$ is the limit in law of $m^a$ for some $\Omega$ and $\xi$, then:
\begin{longlist}[(ii)]
\item[(i)] $\Eb{e^{t m}} < \infty$;
\item[(ii)] furthermore, as $a\to0$,
$\Eb{e^{t m^a}} \to\Eb{e^{t m}}$.
\end{longlist}
\end{corollary}

The proof is straightforward. Note that for any $t\in\R$,
by Fatou's lemma one has that
%
%
\begin{equation}
\label{e.fatou} \E \bigl[e^{t m}\bigr] \le\liminf_{a\to0}
\E \bigl[e^{t m^a}\bigr],
\end{equation}
which implies (i).
Now (ii) follows easily from (i) (used with some $\tilde t > |t|$
and with $\xi=+$), FKG,
and the weak convergence of $m^a$ to $m$.

\subsection{Computing the characteristic function}
Let us prove Proposition~\ref{pr.chf} assuming Proposition~\ref{pr.close} below.
Let $f\in\Sob^3$ be fixed once and for all. For any $k\geq1$, note that
\begin{eqnarray*}
\bigl|\E \bigl[\bigl\langle\Phi^a, f\bigr\rangle^k\bigr]\bigr| &=& \biggl|\E \biggl[
\biggl( \sum_{x\in[0,1]^2 \cap a \Z^2} a^{15/8} f(x)
\sigma_x \biggr)^k\biggr]\biggr|
\\
&=& \biggl| a^{15k/8} \sum_{z_1, \ldots,z_k \in([0,1]^2\cap a\Z^2)^k}
f(z_1) \cdots f(z_k) \Eb{ \sigma_{z_1}
\cdots\sigma_{z_k}} \biggr|
\\
&\le& \|f\|_\infty^k a^{15k/8} \sum
_{z_1, \ldots,z_k \in([0,1]^2\cap a\Z^2)^k} \Eb{ \sigma_{z_1}\cdots
\sigma_{z_k}}
\\
&\le& \| f\|_\infty^k \E \bigl[ \bigl(m^a
\bigr)^k\bigr].
\end{eqnarray*}

Since, by Proposition~\ref{pr.expomoments}, $m^a$ has uniform
exponential moments, we deduce that the series
\[
\phi_{\Phi^a}(f) = \E \bigl[e^{i \langle\Phi^a, f\rangle}\bigr] = 1+ \sum
_{k\geq1} \frac{i^k \Eb{\langle\Phi^a, f\rangle^k}}{k !}
\]
is indeed summable.
Now, for each $k\geq1$, let us prove that the kth moment $\Eb{\<{\Phi
^a, f}^k}$ has a limit as $a\to0$.
Let us fix some cut-off $\varepsilon>0$ and let us divide the kth moment
as follows:
%
%
\begin{eqnarray}
\label{e.split} \qquad \E \bigl[\bigl\langle\Phi^a, f\bigr\rangle^k \bigr]
&=& \sum_{z_1, \ldots,z_k \in([0,1]^2\cap a\Z^2)^k} a^{15k/8} f(z_1)
\cdots f(z_k) \Eb{ \sigma_{z_1} \cdots
\sigma_{z_k}}
\nonumber
\\
&=& \mathop{\sum_{z_1,\ldots, z_k}}_{|z_i - z_j|\geq\eps\ \forall i\neq
j }
a^{2k} f(z_1)\cdots f(z_k) a^{-k/8}
\Eb{ \sigma_{z_1} \cdots\sigma_{z_k}}
\\
&&{}+ \mathop{\sum_{z_1,\ldots, z_k }}_{\inf_{i\neq j }|z_i - z_j| < \eps}
a^{15k/8 } f(z_1)\cdots f(z_k) \Eb{
\sigma_{z_1}\cdots\sigma_{z_k}}.
\nonumber
\end{eqnarray}

Using Theorem~\ref{th.hongler} and assuming Wu's result, we have that
in the domain $[0,1]^2$, there exists a function $z_1,\ldots,z_k
\mapsto\langle z_1,\ldots,z_k\rangle_{[0,1]^2}$ such that
%
%
\begin{equation}
a^{-k/8} \EFK{[0,1]^2_a} {} {
\sigma_{z_1}\cdots\sigma_{z_k}} \longrightarrow \langle z_1,
\ldots,z_k\rangle_{[0,1]^2},
\end{equation}
uniformly in $\inf_{i\neq j} |z_i-z_j| \geq\eps$ [again, up to a
change by a deterministic scalar in the definition of these functions
which arises from normalizing by either $\varrho(a)^{k/2}$ or
$a^{k/8}$]. The fact that the convergence is uniform implies that the
first term in equation~\eqref{e.split} converges as the mesh size $a\to
0$ to
\[
\int\!\!\mathop{\int_{z_1,\ldots, z_k \in([0,1]^2)^k}}_{|z_i - z_j|\geq
\eps\ \forall i\neq j } f(z_1)
\cdots f(z_k) \langle z_1,\ldots,z_k
\rangle_{[0,1]^2} \,dz_1\cdots\, dz_k.
\]

To conclude the proof, it remains to prove that the second term in
equation~\eqref{e.split} is small uniformly in $0<a<\eps$,
when the cut-off $\eps$ is small. This is the content of the next section.

\subsection{Handling the ``local'' $k$-tuples}

%
\begin{proposition}\label{pr.close}
Let $\Omega$ be a domain with $+$ boundary conditions. For any $k\geq
1$, there exist constants
$C_k=C_k(\Omega)< \infty$ such that, for all $0<a<\eps$,
%
\[
\sum_{(x_1,\ldots, x_k)\dvtx \inf_{i\neq j} \{ |x_i - x_j|\} <\eps} a^{15k/8} \E \Biggl[ \prod
_1^k \sigma_{x_i} \Biggr] \le C_k
\eps^{7/4}.
\]
\end{proposition}

\begin{pf}
Our proof is based on the FK representation; we remark that a somewhat
different proof can be obtained by using the Gaussian correlation
inequalities of \cite{\Chuck}.
One implements the $+$ boundary condition via a ghost vertex
corresponding to the boundary and then reduces estimates of $k$th
moments essentially to one
and two point correlations. Those are handled by arguments like in the
Appendix~\ref{s.appendix} below; see especially equation~(\ref{eq.2point}). We now
proceed with more details
using the FK representation approach.

One can write $\Eb{\prod_1^k \sigma_{x_i}}$ using FK as follows:
let $\Delta_k$ be the set of graphs $\Gamma$ defined on the set of
vertices $V_k:= \{1,\ldots,k\} \cup\{ + \}$,
and which are such that the clusters of $\Gamma$ which do not contain
the point $+$ are all of even size.
(Of course, the number $|\Delta_k|$ of such graph structures is finite.)

Now, similarly to Wick's theorem, one has the identity
\[
\E \Biggl[\prod_1^k \sigma_{x_i}\Biggr]
= \sum_{\Gamma\in\Delta_k} \P \bigl[ A_{x_1,\ldots, x_k}(\Gamma)\bigr],
\]
where $A_{x_1,\ldots,x_k}(\Gamma)$ is the event that the graph
structure induced by the FK configuration $\omega$ on the set
$\{x_1,\ldots, x_k\} \cup\{ \partial\Omega\}$ is given by the graph
$\Gamma\in\Delta_k$.


Note that if $\Gamma$ 
is not connected, there is some negative information inherent to the event
$A_{x_1,\ldots, x_k}(\Gamma)$. To overcome this, let $\overline
A_{x_1,\ldots, x_k}(\Gamma)$ be the event that the
graph induced by the FK configuration $\omega$ on $\{x_1,\ldots, x_k\}
\cup\{ \partial\Omega\}$ \textit{includes} the graph $\Gamma$.
Defined this way, $\overline A_{x_1,\ldots, x_k}(\Gamma)$ is an
increasing event (which will allow us to use FKG) and one has for any
$\Gamma\in\Delta_k$:
\[
\P \bigl[A_{x_1,\ldots,x_k}(\Gamma)\bigr] \le\P \bigl[\overline A_{x_1,\ldots,x_k}(\Gamma)\bigr].
\]

Therefore, it is enough for us to prove the following upper bound:
\[
\sum_{(x_1,\ldots, x_k)\dvtx \inf_{i\neq j} \{ |x_i - x_j|\} <\eps} \sum_{\Gamma\in\Delta_k}
\P \bigl[ \overline A_{x_1,\ldots, x_k}(\Gamma)\bigr] \le C \eps^{7/4}
a^{-(15k)/8}.
\]

This is the subject of the next lemma, which concludes the proof of the proposition.
\end{pf}

%
\begin{lemma}
For any domain $\Omega$ and any $k\geq1$, there exists a constant $C_k
= C_k(\Omega) < \infty$ such that,
for all $0<a<\eps$, one has
\begin{longlist}[(ii)]
\item[(i)]
$
\displaystyle \sum_{x_1,\ldots, x_k \in\Omega_a} \sum_{\Gamma\in\Delta_k}
\P\bigl[\overline A_{x_1,\ldots, x_k}(\Gamma)\bigr] \le C_k a^{-(15k)/8},
$\vspace*{10pt}

\item[(ii)]
$
\displaystyle \sum_{(x_1,\ldots, x_k)\dvtx \inf_{i\neq j} \{ |x_i - x_j|\} <\eps} \sum_{\Gamma\in\Delta_k}
\P \bigl[ \overline A_{x_1,\ldots, x_k}(\Gamma)\bigr] \le C_k \eps^{7/4}
a^{-(15k)/8}.
$
\end{longlist}
\end{lemma}

\begin{pf*}{Proof (Sketch)}
The proof of this lemma proceeds by induction. For $k=1$, the bounds
follow easily from Proposition~\ref{pr.1and2} in the Appendix~\ref{s.appendix}.
For $k=2$, using again Proposition~\ref{pr.1and2} and summing $\Pb{x_1
\leftrightarrow x_2}$ over all $x_1,x_2$ which are such
that $|x_1-x_2| \in(2^{-b-1}, 2^{-b}]$, one\vspace*{1pt} gets a bound of the form
$O(1)a^{-4}2^{-2b}(2^{-b}/a)^{-1/4} = O(1) 2^{-7b/4} a^{-15/4}$,
where $a^{-4}2^{-2b}=a^{-2}(2^{-b}/a)^2$ comes from the number of ways
one can choose $x_1$ and $x_2$.
Summing over all possible values of $b$ smaller than $\log_2(a^{-1})$
gives the first bound, while summing
over values of $b$ such that $\log_2(\eps^{-1}) \le b \le\log
_2(a^{-1})$ gives the second bound.
(We neglect boundary issues that can easily be dealt with.)

Let now $k\geq3$ and assume that property (i) holds for all $k'<k$.
We will first
prove that it implies property (ii) from which (i) easily follows
[in fact formally (i)
readily follows from (ii) by taking $\eps$ large enough but due to
boundary issues,
it is better to divide the study into these two sums].


The outer sum in (ii) is over the ordered $k$-tuples $(x_1,\ldots,x_k)$ which are such that
$l:= \inf_{i\neq j} |x_i - x_j| < \eps$.
For any such $k$-tuple $(x_1,\ldots, x_k)$, let us choose one point
among all points which are
at distance $\inf_{i\neq j} |x_i - x_j|$ from at least one of the others
(there are at most $k$ ways to pick one) and let us reorder the points
into a $k$-tuple
$(\hat x_1, \ldots, \hat x_k)$ so that the point we have chosen is
$\hat x_1$.

This way, we obtain
\begin{eqnarray*}
&&\mathop{\sum_{x_1,\ldots, x_k \in\Omega_a}}_{\inf_{i\neq j} |x_i -
x_j| < \eps
} \sum
_{\Gamma\in\Delta_k} \P \bigl[ \overline A_{x_1,\ldots, x_k}(\Gamma)\bigr]
\\
&&\qquad \le k \mathop{\sum_{\hat x_1, \ldots, \hat x_k}}_{\inf_{i\neq j} |\hat x_i
- \hat x_j| = \inf_{i\neq1}|\hat x_1 - \hat x_i| < \eps
} \sum
_{\Gamma\in\Delta_k} \P\bigl[\overline A_{\hat x_1,\ldots, \hat
x_k}(\Gamma)\bigr].
\end{eqnarray*}
Now, for any such $(\hat x_1,\ldots, \hat x_k)$, we split the sum over
$\Gamma\in\Delta_k$ in two parts.
\begin{longlist}[(2)]
\item[(1)] Consider first the sum over graphs $\Gamma$ such that the cluster
of $\hat x_1$ in $\Gamma$ contains
a point $\hat x_{m}$ at distance $<2\varepsilon$ from $\hat x_1$.
Again by reordering (and possibly losing a factor of $k$), one can
assume that $\hat x_m = \hat x_2$.
Now let $A_{\hat x_1, \hat x_2}$ be an annulus which surrounds $\hat
x_1$ and $\hat x_2$ and which is such
that, by RSW, there is probability $c>0$ of the event $S=S(A_{\hat x_1,
\hat x_2})$ that there is an open
path in $A_{\hat x_1, \hat x_2}$ surrounding $\hat x_1$ and $\hat x_2$.

Let $\widehat\Gamma$ be a graph on $\{\hat x_3,\ldots, \hat x_k\}$
obtained from $\Gamma$ in the following way.
If the cluster of $\hat x_1$ and $\hat x_2$ in $\Gamma$ does not
contain other points, let $\widehat\Gamma=\Gamma\setminus\{\hat
x_1,\hat x_2\}$.
Otherwise, first add some connection, if necessary, to make the cluster
of $\hat x_1$ and $\hat x_2$ in $\Gamma$
connected without using $\hat x_1$ and $\hat x_2$ (i.e., all other
vertices are connected by paths that do not pass
through ${\hat x}_1$ and ${\hat x}_2$), and then remove $\hat x_1$ and
$\hat x_2$ from the cluster.
Note that, in both cases, $\widehat\Gamma\in\Delta_{k-2}$.
Using FKG, one can easily check that
\begin{eqnarray*}
\P \bigl[\overline A_{\hat x_1,\ldots, \hat x_k}(\Gamma)\bigr] & \le& (1/c) \P\bigl[ \overline
A_{\hat x_1,\ldots, \hat x_k}(\Gamma)\mbox{ and }S\bigr]
\\
&&\qquad \le(1/c) \P \bigl[ \hat x_1 \leftrightarrow\hat x_2
\mbox{ and }\overline A_{\hat x_3,\ldots, \hat x_k}(\widehat\Gamma )\mbox{ and } S\bigr]
\\
&&\qquad \le(1/c) \FK{} {+} {\hat x_1 \leftrightarrow\hat
x_2 } \P\bigl[ \overline A_{\hat x_3,\ldots, \hat x_k}(\widehat\Gamma) \bigr]
\\
&&\qquad \le O(1) \,d^{-1/4} a^{1/4} \P \bigl[\overline
A_{\hat x_3,\ldots, \hat x_k}( \widehat\Gamma) \bigr],
\end{eqnarray*}
where $d$ denotes the distance between $\hat x_1$ and $\hat x_2$ and by
$+$ we mean 
wired b.c. on the inner boundary of $A_{\hat x_1, \hat x_2}$.

Summing over all $\hat x_1,\ldots, \hat x_k$ which are such that
$d=|\hat x_1 - \hat x_2| \in(2^{-b-1}, 2^{-b}]$,
and considering that there are at most $k^2$ ways of choosing $\hat
x_1$ and $\hat x_2$ from $\{x_1, \ldots, x_k\}$,
this case gives a contribution which is bounded by
\[
O(1) k^2 a^{-4} 2^{-2b} 2^{b/4}
a^{1/4} C_{k-2} a^{-(15 (k-2))/8},
\]
where $a^{-2}(2^{-b}/a)^2=a^{-4} 2^{-2b}$ is an upper bound on the
number of ways to choose $\hat x_1$ and $ \hat x_2$ from $\Omega_a$.
Hence, we get the following upper bound:
\[
O(1) k^2 2^{-7b/4} C_{k-2} a^{-(15k)/8}.
\]
It remains to sum over the possible values of $b$, that is, $\log_2(\eps
^{-1}) \le b \le\log_2(a^{-1})$, which gives a bound of the desired form.

Note that we neglected boundary issues here (they can be handled easily
at least if $\partial\Omega$ is smooth enough).

\item[(2)] Consider now the remaining sum over graphs $\Gamma$ such that the
cluster of $\hat x_1$ in $\Gamma$ does not contain
any point at distance $<2\varepsilon$ from $\hat x_1$. In this case, there
is at least one point, say $\hat x_2$, which is at distance
$l$ from $\hat x_1$. If the cluster of $\hat x_2$ in $\Gamma$ contains
a point at distance $<2\varepsilon$ from $\hat x_2$, then we
can take $\hat x_2$ to play the role of $\hat x_1$ and we are back in
situation 1.
We can therefore assume that the cluster of $\hat x_2$ in $\Gamma$ does
not contain any point at distance $<2\varepsilon$ from
$\hat x_2$. We can then pick an annulus $A_{\hat x_1,\hat x_2}$ that
surrounds $\hat x_1$ and $\hat x_2$ and does not contain
any other point belonging to the clusters of $\hat x_1$ and $\hat x_2$
in $\Gamma$, and which, by RSW, contains an open
path surrounding $\hat x_1$ and $\hat x_2$ with probability $c>0$. We
call $S=S(A_{\hat x_1,\hat x_2})$ the latter event.
If $S$ occurs, $\hat x_1$ and $\hat x_2$ belong to the same FK cluster.
If we denote by $\widehat\Gamma$ a graph on
$\{\hat x_3,\ldots, \hat x_k\}$ obtained from $\Gamma$ by connecting
the clusters of $\hat x_1$ and $\hat x_2$ in $\Gamma$
outside of $\hat x_1$ and $\hat x_2$, and then removing $\hat x_1$ and
$\hat x_2$ from $\Gamma$, we have that
$\widehat\Gamma\in\Delta_{k-2}$. Using FKG, one can easily check that
\begin{eqnarray*}
\P \bigl[\overline A_{\hat x_1,\ldots, \hat x_k}(\Gamma)\bigr] & \le& (1/c) \P \bigl[\overline
A_{\hat x_1,\ldots, \hat x_k}(\Gamma)\mbox{ and } S\bigr]
\\
& \le& (1/c) \P \bigl[ \hat x_1 \leftrightarrow\hat x_2
\mbox{ and }\overline A_{\hat x_3,\ldots, \hat x_k}(\widehat\Gamma )\mbox{ and } S\bigr]
\\
& \le& (1/c) \FK{} {+} {\hat x_1 \leftrightarrow\hat x_2
} \P \bigl[\overline A_{\hat x_3,\ldots, \hat x_k}(\widehat\Gamma) \bigr]
\\
& \le& O(1) l^{-1/4} a^{1/4} \P\bigl[\overline A_{\hat x_3,\ldots, \hat x_k}(
\widehat\Gamma) \bigr],
\end{eqnarray*}
where by $+$ we mean 
wired on the inner boundary of $A_{\hat x_1, \hat x_2}$.

Summing over all $x_1,\ldots, x_k$ which are such that $l= |\hat x_1 -
\hat x_2| = \inf_{i\neq j} |x_i - x_j| \in(2^{-b-1}, 2^{-b}]$, this
case gives a contribution which is bounded by
\[
O(1) k^2 a^{-4} 2^{-2b} 2^{b/4}
a^{1/4} C_{k-2} a^{-(15 (k-2))/ 8},
\]
where $k^2$ comes from the ways of choosing $\hat x_1$ and $\hat x_2$
from $\{x_1, \ldots, x_k\}$ and $a^{-2}(2^{-b}/a)^2=a^{-4} 2^{-2b}$ is
an upper bound on the
number of ways to choose $\hat x_1$ and $ \hat x_2$ from $\Omega_a$.
Hence, we get the following upper bound:
\[
O(1) k^2 2^{-7b/4} C_{k-2} a^{-(15 k)/ 8}.
\]
%
It remains to sum over the possible values of $b$, that is, $\log_2(\eps
^{-1}) \le b \le\log_2(a^{-1})$, which gives the desired result.

Modulo boundary issues that are easily dealt with, this concludes the
proof of the lemma, which in turn implies the proposition.\quad\qed
\end{longlist}\noqed
\end{pf*}

\subsection{Consequences of this approach}
This proof of Theorem~\ref{th.main} through the study of the moments of
$m^a$ sheds some light on $\Phi^\infty$.
For example, it enables in some cases to explicitly compute the
variance of $m_\infty$. Indeed, in the full plane $\C$,
if one looks at $\langle\Phi, 1_A \rangle$, then from the work of \cite
{arXiv1202.2838} or \cite{arXiv1112.4399}, we get that
%
%
\begin{eqnarray}
\E \bigl[\langle\Phi_\C, 1_A\rangle^2\bigr] &=& C
\int\!\!\int_A \frac{1}{|x-y|^{1/4}} \,dx \,dy.
\end{eqnarray}
Here, $C$ is a constant which can be computed explicitly thanks to the formula
(see Theorem~\ref{th.Wu}) by Wu (\cite{\Wu}).
Therefore, the second moment of $\langle\Phi_\C, 1_A \rangle$ can be
computed numerically or exactly depending on the set $A$.

\section{Without assuming Wu's result}\label{Theta}

The purpose of this section is to briefly explain how to adapt our
proofs if one does not
want to rely on Wu's result, Theorem~\ref{th.Wu}.
In this case, as explained in Section~\ref{s.renormalization}, one
would need to renormalize our fields by
%
%
\begin{equation}
  \Theta_a:= a^2 \varrho(a )^{-1/2},
\end{equation}
instead of $\Theta_a = a^{15/8}$.

\subsection{Adapting the first proof (Section~\protect\ref{first})}
Let us point out here that it is not a priori needed to have an exact
rescaling of the form $a^{15/8}$ if one wants to obtain
our main result, Theorem~\ref{th.main}. For example, this situation
arises in \cite{arXiv1008.1378}, where the four-arm event is only known up to
possible logarithmic corrections. Therefore, in order to build the
\textit{pivotal measures} there, it is not possible to assume a renormalization
of the discrete counting measure by $\eta^{3/4}$; instead, a more
cumbersome renormalization of $\eta^2 \alpha_4(\eta,1)^{-1}$ is
needed---see \cite{arXiv1008.1378} for more details.
In the present work, the same technology as in \cite{arXiv1008.1378} would
enable us to prove Theorem~\ref{th.main} without relying on Wu's result.

Yet, some of the present proofs would need to be slightly modified and
some quantitative lemmas (such as Lemma~\ref{l.333}, e.g.) would
need to be changed.
Let us point out that we would have at our disposal the following
useful bound on the one-arm event:
%
%
\begin{equation}
\label{e.onearmbound} C n^{-1/2} \le\alpha_1^\mathrm{FK}(n)
\le n^{-\alpha}
\end{equation}
for some exponent $\alpha>0$.
The lower bound follows from Smirnov's observable (see \cite{arXiv0912.4253})
while the upper bound follows from the RSW theorem in \cite{arXiv0912.4253}.
Such bounds are enough to carry the proof from \cite{arXiv1008.1378} through
(except for the conformal covariance property, Theorem~\ref{th.CC},
which needs at least an $\mathrm{SLE}_{16/3, 16/3-6}$ computation for the
one-arm event).

\subsection{Adapting the second proof (Section~\texorpdfstring{\protect\ref{second}}{3})}

The second proof is easier to adapt, since the results in \cite{arXiv1202.2838} are stated precisely with the renormalization factor
$\Theta_a = a^2 \varrho(a )^{-1/2}$. The proof of Proposition~\ref{pr.expomoments} works as before.
Of course, Proposition~\ref{pr.close} would be stated in a less
quantitative manner, but using, for example, the above estimate~\eqref
{e.onearmbound}, one could still handle the local $k$-tuples, which
would give us the desired result.

\section{Properties of the limiting magnetization field 
\texorpdfstring{$\Phi^\infty$}{Phiinfty}}\label{properties}

In this last section, we wish to list some interesting properties
satisfied by the magnetization field $\Phi^\infty$
which will be proved in \cite{\CGNproperties} as well as some results
on the near-critical behavior of the Ising
model along the $h$-direction which appear or will appear in \cite
{\CGNNC,\CGNproperties}.

\begin{enumerate}
\item
If $m^\infty_\Omega$ denotes the scaling limit of the renormalized magnetization
in a bounded domain $\Omega$, then there exists a constant $c=c_\Omega
>0$ such that
%
%
\begin{equation}
\label{eq:tail} \log\P \bigl[m^\infty_\Omega> x\bigr] \underset{x\to
\infty} {\sim} -c x^{16}.
\end{equation}
Furthermore, one can show that the constant $c=c_\Omega$ depends on the
domain but does not depend on the boundary conditions.
We point out that \eqref{eq:tail} clearly shows the non-Gaussianity of
the magnetization field; the latter already follows
from the fact that the correlation functions computed by Chelkak,
Hongler and Izyurov in \cite{arXiv1202.2838} do not satisfy
Wick's formula.

\item The probability density function of $m^\infty_\Omega$ is smooth
as a consequence of the following quantitative bound on its Fourier
transform: $\forall t\in\R$,
%
%
\begin{equation}
  \bigl|\E_{\Omega} \bigl[e^{i t m^\infty}\bigr]\bigr| \le e^{-C |t|^{16/15}}
\end{equation}
for some 
constant $C>0$.

\item In \cite{\CGNproperties}, it will be shown that the Ising model
on the rescaled lattice $a \Z^2$ with renormalized external magnetic
field $h_a:= h a^{15/8}$ has a \textit{near-critical} (or \textit
{off-critical}) scaling limit as $a \searrow0$. This near-critical
limit is no longer scale-invariant but is conformally covariant instead
and has exponential
decay of its correlations.

\item Finally, in \cite{\CGNNC} we prove that the average magnetization
$\<{\sigma_0}$ of the Ising model on $\Z^2$ at $\beta=\beta_c$ and with
external magnetic
field $h>0$ is such that
%
%
\begin{equation}
\langle\sigma_0\rangle_{\beta_c, h} \asymp h^{1/15 }.
\end{equation}
\end{enumerate}

\begin{appendix}
\section{Tightness of the magnetization field}\label{tight}
We first introduce the setup for proving tightness
when the field $\Phi^a$ is defined on the compact square $[0,1]^2$.
The extension to general domains as well as to the \textit{full} plane
will be given in Section~\ref{s.extension}.

\subsection{Tightness for \texorpdfstring{$\Phi^a$}{Phia} in a well-chosen 
Sobolev space (case
of the square)}\label{s.tightness}

In this subsection, we follow (almost word for word) the functional
approach which was used by Julien Dub\'edat in \cite{\Dubedat} for
another well-known field: the Gaussian Free Field.

Let $\Sob^1_0 = \Sob_0^1([0,1]^2)$ be the classical Sobolev Hilbert
space, that is, the closure of $C_0^\infty([0,1]^2)$ for the norm
\[
\| f \|_{\Sob^1}^2:= \int_{[0,1]^2} \|\nabla
f \|^2 \,dA.
\]

Let $\Sob^{-1}$ be the dual space of $\Sob^1_0$. It is a space of
distributions (i.e., $\Sob^{-1} \subset D'$) and it is also a Hilbert
space equipped with the norm (the \textit{operator norm} on $\Sob^{-1}$)
\[
\| h \|_{\Sob^{-1}}:= \sup_{g \in C_0^\infty([0,1]^2)\dvtx \| g\|_{\Sob^1}
\le1} \langle h, g \rangle.
\]
(Here, $\langle h, g \rangle$ stands for the evaluation of the
distribution $h$ against the test function~$g$.)

It will be useful to work with the following basis of the space
$C_0^\infty([0,1]^2)$ endowed with the $L^2$ norm:
for any $j,k \in\N^+$, let
%
%
\begin{equation}
e_{j,k}(x,y):= 2 \sin(j \pi x) \sin(k \pi y).
\end{equation}

It is straightforward to check that
%
%
\begin{equation}
\cases{
\displaystyle (e_{j,k})_{j,k>0}\qquad\mbox{is a joint orthogonal basis for } \Sob^{-1}\mbox{ and }\Sob^1_0,
\vspace*{5pt}\cr
\displaystyle \| e_{j,k} \|_{\Sob^1}^2 = j^2 + k^2,
\vspace*{3pt}\cr
\displaystyle \| e_{j,k} \|_{\Sob^{-1}}^2 = \frac{1} {j^2 + k^2}.}
\end{equation}
In particular, if $h = \sum_{j,k} a_{j,k} e_{j,k}$, then $\| h\|_{\Sob
^{-1}}^2 = \sum_{j,k} \frac{a_{j,k} ^2} {j^2 + k^2}$.

More generally, for any $\alpha>0$, one can define the Hilbert space
$\Sob^{\alpha}_0$ as the closure of $C_0^\infty([0,1]^2)$ for the norm
\[
\| f \|_{\Sob^\alpha}^2:= \sum_{j,k>0}
a_{j,k}^{2} \bigl(j^2 + k^2
\bigr)^{\alpha},
\]
where $f \in C_0^\infty$ is decomposed as $f = \sum_{j,k>0} a_{j,k} e_{j,k}$.
Let $\Sob^{-\alpha}$ be its dual space. It is a Hilbert space with norm
\begin{eqnarray*}
\| h \|_{\Sob^{-\alpha}} &:=& \sup_{g \in C_0^\infty([0,1]^2)\dvtx \|
g\|
_{\Sob^\alpha} \le1} \langle h, g \rangle.
\end{eqnarray*}

Furthermore, if $h\in L^2 \subset\Sob^{-\alpha}$, then $h$ has a
Fourier expansion and its $\|\cdot\|_{\Sob^{-\alpha}}$ norm can be
expressed as
%
%
\begin{equation}
\label{e.Sob} \| h \|_{\Sob^{-\alpha}}^2 = \sum
_{j,k} a_{j,k}^2 \frac{1} {(j^2 + k^2)^\alpha}.
\end{equation}


We will also make use of the following classical result.

\begin{proposition}[(Rellich theorem)]
For any $\alpha_1<\alpha_2$, $\Sob^{-\alpha_1}$ is compactly embedded
in $\Sob^{-\alpha_2}$ ($\Sob^{-\alpha_1} \subset\subset\Sob^{-\alpha_2}$).
In particular, for any $R>0$,
the ball
\[
\overline{B_{\mathcal{H}^{-2}}(0,R)}
\]
is compact in $\mathcal{H}^{-3}$.
\end{proposition}

Thanks to this property, in order to prove tightness, it is enough for
us to prove the following result.

\begin{proposition}\label{pr.tight}
Let us fix some boundary condition $\xi$ on the square $[0,1]^2$.
Assume that the boundary condition $\xi$
is made of finitely many
arcs of~$+,-$ or \textit{free} type. By $\Phi^a$, we denote the
magnetization field within $[0,1]^2 \cap a \Z^2$
subject to the boundary condition $\xi$.Then as $a\searrow0$, one has
\[
\limsup_{a\searrow0} \E\bigl[ \bigl\| \Phi^a
\bigr\|_{\Sob^{-2}}^2 \bigr] < \infty,
\]
uniformly in the boundary conditions $\xi$, and thus
$(\Phi^a)_{a>0}$ is \textit{tight} in the space~$\Sob^{-3}$.
\end{proposition}



\begin{pf}
We wish to bound from above the quantity
\begin{eqnarray*}
\E \bigl[ \bigl\| \Phi^a \bigr\|_{\Sob^{-2}}^2 \bigr] & =& \E \biggl[ \sum
_{j,k >0} \bigl\langle\Phi^a,
e_{j,k} \bigr\rangle^2 \frac{1}{(j^2 +k^2)^2}\biggr]
\\
&=& \sum_{j,k>0} \frac{1} {(j^2 +k^2)^2} \E \bigl[ \bigl\langle
\Phi^a, e_{j,k} \bigr\rangle^2\bigr].
\end{eqnarray*}
This is clone using the following lemma.
\end{pf}

%
\begin{lemma}\label{l.SOB}
There is a constant $C>0$ such that for all $j,k>0$
\[
\limsup_{a\to0} \sup_{j,k} \E \bigl[ \bigl\langle
\Phi^a, e_{j,k} \bigr\rangle^2\bigr] < C.
\]
\end{lemma}

\begin{pf}
\begin{eqnarray*}
&& \E \bigl[ \bigl\langle\Phi^a, e_{j,k} \bigr
\rangle^2\bigr]
\\
&&\qquad \le  a^{15/4} \sum_{x \neq y \in[0,1]^2 \cap a \Z^2}
\biggl| \int\!\!\int_{S_a(x) \times S_a(y)} \frac{\Eb{\sigma_x \sigma_y}} {a^4} e_{j,k}(\bar x)
e_{j,k}(\bar y) \,dA(\bar x) \,dA(\bar y) \biggr|
\\
&&\quad\qquad {}+ a^{15/4} \sum_{x\in[0,1]^2 \cap a \Z^2} \biggl( \int
_{S_a(x)} \frac{1} {a^2} e_{j,k}(x) \,dA(\bar x)
\biggr)^2.
\\
&&\qquad \le a^{15/4} \| e_{j,k}
\|_\infty^2 \sum_{x\neq y \in[0,1]^2 \cap a \Z^2} \bigl|\Eb{
\sigma_x \sigma_y}\bigr| + a^{15/4} \|
e_{j,k} \|_\infty^2 \sum
_{x \in[0,1]^2 \cap a \Z^2} 1
\\
&&\qquad \le O(1),
\end{eqnarray*}
uniformly in $j,k$ and the boundary condition $\xi$ (indeed, by FKG for
the FK representation,
it is enough to dominate $\Eb{\sigma_x \sigma_y}$ by the extreme
boundary conditions $\xi=+$ or $\xi=-$).
\end{pf}


%
\begin{remark}\label{r.baralpha}
Using Lemma~\ref{l.SOB} as is, it is straightforward to strengthen the
above proposition by showing\vspace*{1pt} that for any $\eps>0$,
$(\Phi^a)_{a>0}$ is in fact \textit{tight} in the space $\Sob^{-1-\eps}$.
It is thus natural to wonder for which values of $\alpha>0$, $(\Phi
^a)_{a>0}$ remains tight in $\Sob^{-\alpha}$.
It is clear that there is a lot of room if one wishes to obtain better
estimates than the one provided by Lemma~\ref{l.SOB}.
Yet it appears that there is some $\bar\alpha>0$ such that $(\Phi
^a)_{a>0}$ is not tight in $\Sob^{-\alpha}$ when $\alpha<\bar\alpha$.
In particular, it appears that
$\Phi^{\infty}=\lim_{a \to0} \Phi^a$ is less regular than the planar
Gaussian free field.
\end{remark}

\subsection{Extension to other domains and to the full plane}\label{s.extension}





The next subsection is concerned with the case of bounded domains;
later we will tackle the case of the infinite plane.

\subsubsection{Case where \texorpdfstring{$\Omega\subsetneqq\C$}{$OmegasubsetneqqC$} is a bounded simply
connected domain of the plane, with
prescribed boundary condition \texorpdfstring{$\xi$}{xi} on
\texorpdfstring{$\partial\Omega$}{partialOmega}}

Let $(Q_i)_{i\in\N}$ be a \textit{Whitney decomposition} of $\Omega$ into
disjoint squares.
For any $a>0$, let $\Phi_\Omega^a$ be the magnetization field on $\Omega
\cap a \Z^2$ induced by the boundary condition $\xi$.
One can write $\Phi^a$ as
\[
\Phi^a = \sum_{i\in\N}
\Phi^a_{\md Q_i}.
\]
By the triangle inequality, one has that
\[
\bigl\| \Phi^a \bigr\|_{\Sob^{-2}} \le\sum_{i\in\N}
\bigl\| \Phi^a_{\md Q_i} \bigr\|_{\Sob^{-2}}.
\]

Now the key step is to notice that the proof of Proposition~\ref{pr.tight} immediately gives the following.

\begin{lemma}\label{l.333}
There exists a uniform constant $C>0$ such that for any domain~$\Omega$
and any boundary condition $\xi$ on $\partial\Omega$,
if $Q_i$ is a square inside $\Omega$ [with area $\lambda(Q_i)$], then
for any $a>0$:
%
%
\begin{equation}
\E \bigl[\bigl\| \Phi^a_{\md Q_i}\bigr\|^2_{\Sob^{-2}}\bigr] \le C
\lambda(Q_i)^{15/8}.
\end{equation}
\end{lemma}

\begin{pf*}{Proof (Sketch)}
To see why this holds, take a square $Q_i$ inside $\Omega$. Let $q$ be
its side-length so that $q^2 = \lambda(Q_i)$.
By renormalizing the scale by a factor $1/q$, one can see that our
field $\Phi^a_{\md Q_i}$ has the same $\Sob^{-2}$ norm as
\[
q^{15/8} \times\Phi^{a/q}_{\md(1/q) Q_i}.
\]
But now, $\frac{1} q Q_i$ is a square of side-length 1, therefore by
Proposition~\ref{pr.tight} (which was uniform in the outer boundary condition)
\[
\E \bigl[\bigl\| \Phi^{a/q}_{\md(1/q) Q_i} \bigr\|_{\Sob^{-2}}^2\bigr] \le
C.
\]
This gives
\[
\E \bigl[\bigl\| \Phi^{a}_{\md Q_i} \bigr\|_{\Sob^{-2}}^2\bigr]
\le q^{15/4} C = C \lambda(Q_i)^{15/8}.
\]\upqed
\end{pf*}

By Cauchy--Schwarz, this implies that
%
%
\begin{equation}
\label{e.weakarea} \E \bigl[\bigl\| \Phi^a_{\md Q_i}\bigr\|_{\Sob^{-2}}\bigr] \le
C^{1/2} \lambda(Q_i)^{15/16}.
\end{equation}
From this formula, one can see that one cannot hope to prove a
tightness result for $\Phi^a$ on the full domain $\Omega$.
Indeed there are bounded domains for which
$\sum_i \lambda(Q_i)^{15/16}$ diverges.
Yet, for our purposes, it will be sufficient to prove the following
weaker result.

%
\begin{proposition}\label{pr.tight2}
Let $\Omega$ be a bounded simply connected domain of the plane.
For any open set $U$ whose closure $\overline{U}$ is contained inside
$\Omega$, there is a constant
$C=C_U>0$ such that for any boundary condition $\xi$ on $\partial\Omega
$, one has
\[
\E \bigl[\bigl\| \Phi^a_{\md U} \bigr\|_{\Sob^{-2}}\bigr] <
C_U.
\]
Hence, the restriction of $(\Phi^a)_{a>0}$ to the open subset $U$ is a
tight sequence in $\Sob^{-3}$.
\end{proposition}

\begin{pf}
Observe that
\begin{eqnarray*}
\E \bigl[\bigl\| \Phi^a_{\md U} \bigr\|_{\Sob^{-2}}\bigr] & \le&
C^{1/2} \sum_{i, Q_i \cap U \neq\varnothing} \lambda(Q_i)^{15/16}.
\end{eqnarray*}
By the properties of Whitney decompositions, only finitely many $Q_i$
intersect the subset $U$, hence the above sum is finite and is bounded
from above by some constant $C=C(U)>0$.
\end{pf}

\subsubsection{Case of the infinite plane}
(The case of nonbounded simply connected domains is treated similarly.)


Our magnetization field $\Phi^a:= \sum_{x\in a \Z^2} a^{15/8} \sigma_x
\delta_x$ is well defined as a distribution on the full plane $\R^2$.
One natural way to proceed in order to keep some tightness is to view
our field as a nested sequence of restricted fields: $(\Phi^a_{\md
B_k})_{k\geq1}$
where $B_k$ is the square $[-2^k, 2^k]^2$. This sequence of nested
distributions lives in the product of Hilbert spaces
\[
\Sob^{-3}_\infty:= \prod_{k \geq1}
\Sob^{-3}_{B_k},
\]
where for each $k\geq1$, $\Sob^{-3}_{B_k}$ denotes the dual of $\Sob^3_0(B_k)$.

Since for any $k\geq1$, $(\Eb{ \| \Phi^a_{\md B_k} \|_{\Sob
^{-2}_{B_k}} })_{a>0}$ is a bounded sequence [by\break $O(2^{15k/8})$],
the sequence of random variables $(\Phi^a_{\md B_k})_{a>0}$ is tight in
the space $\Sob^{-3}_{B_k}$. In particular, there is a subsequential
scaling limit,
that is, there is a random field $\Phi_k\in\Sob^{-3}_{B_k}$ and a
sequence $(a_m^k)_{m\geq1}$ with $a_m^k \searrow0$ such that
\[
\Phi^{a_m^k}_{\md B_k} \underset{m\to\infty} {\stackrel{d} {
\longrightarrow}} \Phi_k,
\]
in law (for the topology on $\Sob^{-3}_{B_k}$ induced by $\| \cdot\|
_{\Sob^{-3}}$).
Furthermore, from $k$ to $k+1$, one can choose the subsequential
scaling limit
$(a^{k+1}_m)_{m\geq1}$ so that $\{ a^{k+1}_m \}_m \subset\{a^k_m\}_m$.
This allows us to define a ``joint'' subsequential scaling limit along
the sequence
\[
\bar a_m:= a^m_m.
\]
Doing so, the sequence $(\Phi^{\bar a_m}_{\md B_k})_{k\geq1}$
converges in law (for the product topology) to
\[
(\Phi_k)_{k\geq1} \in\Sob^{-3}_\infty.
\]

It is obvious (going back to the discrete mesh fields $\Phi^a_{\md
B_k}$) that a.s. for any $k\geq1$, one has
\[
\Phi_{k+1} 1_{B_k} \equiv\Phi_k.
\]

\section{First and second moments for~the~magnetization}\label{s.appendix}
The main purpose of this appendix is to prove the following proposition
on the first and second moments of the magnetization in a bounded
smooth domain~$\Omega$. (In fact, to simplify the notation, we will
only prove it in the case where $\Omega$ is a square domain; see
Proposition~\ref{pr.1and22}.) Along the way, we will also prove some
useful bounds on the one-arm event in critical FK percolation (Lemma
\ref{l.onearm}).
Let us point out that in this appendix, we do not need to assume
Theorem~\ref{th.Wu}. 

%
\begin{proposition}\label{pr.1and2}
Let $\Omega$ be a bounded smooth domain of the plane.
Let $M^a_\Omega=M^a$ be the (nonrenormalized) magnetization
\[
M^a= \sum_{x\in\Omega_a} \sigma_x.
\]

There is a 
constant $C>0$ such that for each mesh size $a>0$, one has
\begin{longlist}[(ii)]
\item[(i)] $\EFK{}{+}{M_a} \le C a^{-2} \sqrt{\varrho(a)}$
and\vspace*{1pt}

\item[(ii)] $\EFK{}{+}{(M_a)^2} \le C a^{-4} \varrho(a) $.
\end{longlist}
\end{proposition}

[Obviously here (i) follows from (ii) using Cauchy--Schwarz.]
For simplicity of presentation, we will prove this result only in the
particular case where $\Omega$ is a square domain.
Furthermore, in order to simplify the notation in the proof, we will
work with a nonrenormalized lattice.
Before restating the above proposition in this setting, let us
introduce the following notation: for any $N\geq1$, let
%
%
\begin{equation}
  \rho(N):= \EFK{\Z^2} {} {\sigma_{(0,0)}
\sigma_{(N,N)}}.
\end{equation}
As such, $\rho(N)$ is related to $\varrho( a=\sqrt{2} (N)^{-1})$, where
$\varrho(a)$ was defined in~\eqref{e.varrho}.

We will show the following proposition.

%
\begin{proposition}\label{pr.1and22}
For any $N\geq1$, let $\Lambda_N$ be the square $[-N,N]^2$ and let
$M_N$ be the magnetization in $\Lambda_N$, that is,
\[
M_N:= \sum_{x\in\Lambda_N} \sigma_x.
\]
Then there is a 
constant $C>0$ such that for all $N\geq1$,
\begin{longlist}[(ii)]
\item[(i)] $\EFK{}{+}{M_N} \le C N^2 \rho(N)^{1/2}$
and\vspace*{1pt}

\item[(ii)] $ \EFK{}{+}{M_N^2} \le C N^4 \rho(N)$.
\end{longlist}
\end{proposition}

The proof of the proposition relies on the following lemma, which
already appeared in \cite{\CamiaNewmanIsing}.\vadjust{\goodbreak}
To be self-contained, we include a proof here. (Also, the lemma below
includes more than what is actually needed for Proposition~\ref{pr.1and22} but it will be useful for future reference.) We denote by
$\FK{p_c}{\mathrm{free}}{\cdot}$ (resp., $\FK{p_c}{+}{\cdot}$) the
critical FK percolation measure with free (resp., wired) boundary conditions.

%
\begin{lemma}\label{l.onearm}
There exists a constant $C< \infty$ such that
\[
\cases{
\displaystyle \frac{1} C \sqrt{\rho(N)} \le\FK{p_c} {+} {0 \leftrightarrow\partial\Lambda_N} \le C \sqrt{\rho(N)},
\vspace*{5pt}\cr
\displaystyle \frac{1} C \sqrt{\rho(N)} \le\FK{p_c} {\mathrm{free}} {0\leftrightarrow\partial\Lambda_N} \le C \sqrt{\rho(N)},
\vspace*{5pt}\cr
\displaystyle\rho(N)\le C \rho(2N).}
\]
\end{lemma}

\begin{pf}
To derive the first two parts of the lemma, it is clearly enough to
prove the following inequality for some constant $C< \infty$:
\[
  \frac{1} C \sqrt{\rho(N)} \le\FK{p_c} {
\mathrm{free}} {0 \leftrightarrow\partial\Lambda_N} \le
\FK{p_c} {+} {0 \leftrightarrow\partial\Lambda_N} \le C
\sqrt{\rho(N)}.
\]
Let us first handle the LHS: clearly, using FKG, one has
\begin{eqnarray*}
  \rho(N) & \leq&\FK{p_c} {+} {0 \leftrightarrow\partial
\Lambda_{N/2}}^2.
\end{eqnarray*}
Now we wish to show that
%
%
\begin{equation}
\label{e.NN2} \FK{p_c} {+} {0 \leftrightarrow\partial
\Lambda_{N/2}} \le c \FK{p_c} {\mathrm{free} } {0
\leftrightarrow\partial\Lambda_{N}}
\end{equation}
for some constant $c< \infty$.
This can be seen as follows: let $R_N$ be the event that there is open
circuit in the annulus $\Lambda_{N/2}\setminus\Lambda_{N/4}$, then
\begin{eqnarray*}
  \FK{p_c} {\mathrm{free} } {0 \leftrightarrow\partial
\Lambda_{N}} & \geq&\FK{p_c} {\mathrm{free} } {0
\leftrightarrow\partial\Lambda_{N}; R_N}
\\
&\geq&\FK{p_c} {\mathrm{free} } {0 \leftrightarrow\partial
\Lambda_{N} | R_N}
\\
&\geq&\FK{p_c} {+} {0 \leftrightarrow\partial\Lambda_{N/2}}
\FK{p_c} {+ } {\partial\Lambda_{N/4} \leftrightarrow
\partial\Lambda_{N}},
\end{eqnarray*}
which concludes the proof of~\eqref{e.NN2} by using RSW from \cite
{arXiv0912.4253}. Altogether this proves the LHS
inequalities in the first two parts of Lemma~\ref{l.onearm}.
The RHS is proved along the same lines. Namely, one clearly has by FKG that
\begin{eqnarray*}
  \rho(N) & \geq&\FK{p_c} {\mathrm{free}} {0
\leftrightarrow\partial\Lambda_{N/2}}^2.
\end{eqnarray*}
Now obviously, $\FK{p_c}{\mathrm{free}}{0 \leftrightarrow\partial
\Lambda_{N/2}} \geq\FK{p_c}{\mathrm{free}}{0 \leftrightarrow\partial
\Lambda_{2N}}$
and thus, using again\break \eqref{e.NN2}, this concludes the proof of the
first two parts of Lemma~\ref{l.onearm}.
It is easy to see from the above computation that, possibly by changing
the value of $C$, one can get the last
part of Lemma~\ref{l.onearm}.
\end{pf}

\begin{pf*}{Proof of Proposition~\ref{pr.1and22}}
Even though, as pointed out above, property~(i) follows from property~(ii) by Cauchy--Schwarz, we will give a detailed proof of (i) and only
briefly highlight how to deal with (ii).

We divide the domain $\Lambda_N$ into $n\asymp\log_2 N$ disjoint
annuli $A_0, \ldots, A_n$ such that for each $i\in[0,n]$, the vertices\vadjust{\goodbreak}
in $A_i$ are at distance $2^i$ (up to a factor of 2) from
the boundary $\partial\Lambda_N$. This decomposition gives us
\begin{eqnarray*}
  \EFK{} {+} {M_N} & =& \sum
_{0\le i \le n} \sum_{x\in A_i} \FK{} {+} {x
\leftrightarrow\partial\Lambda_N}
\\
&\le& O(1) \sum_{0\le i \le n} \#\{A_i \} \FK{}
{+} {0\leftrightarrow\partial\Lambda_{2^i}}
\\
&\le& O(1) \sum_{0\le i \le n} N 2^i \FK{} {+}
{0\leftrightarrow\partial\Lambda_{2^i}}.
\end{eqnarray*}

Now, one has that for any $i\leq n$,
\begin{eqnarray*}
  \FK{} {+} {0\leftrightarrow\partial\Lambda_{N}} &
\geq&\FK{} {\mathrm{free}} {0\leftrightarrow\partial\Lambda_{2^i}} \FK{}
{\mathrm{free}} {\partial\Lambda_{2^i} \leftrightarrow\partial
\Lambda_{N}}
\\
&\geq& 1/C^2 \FK{} {+} {0\leftrightarrow\partial
\Lambda_{2^i}} \FK{} {\mathrm{free}} {\partial\Lambda_{2^i}
\leftrightarrow\partial\Lambda_{N}}
\end{eqnarray*}
from Lemma~\ref{l.onearm}. Continuing the above computation, one obtains
\begin{eqnarray*}
  \EFK{} {+} {M_N} &\le& O(N) \sum
_{0\le i \le n} C^3 2^i \frac{\sqrt{\rho(N)}} {\FK{}{\mathrm
{free}}{\partial\Lambda_{2^i} \leftrightarrow\partial\Lambda_{N}}}.
\end{eqnarray*}
It is known from \cite{arXiv0912.4253}, Proposition 24, that $\FK{}{\mathrm
{free}}{\partial\Lambda_{2^i} \leftrightarrow\partial\Lambda_{N}}
\geq c (2^{i}/N)^{1/2}$ for some constant $c>0$. This gives
\[
  \EFK{} {+} {M_N} \le O(1) N \sum
_{0\le i \le n} 2^{i/2} N^{1/2} \sqrt{\rho(N)}
\le O(1) N^2 \sqrt{\rho(N)},
\]
which completes the proof of condition (i).

The proof for the second moment (ii) follows exactly the same lines
except that the combinatorics is slightly more tedious.
As an indication, let us give two upper bounds which are useful to
carry out the computation properly: if $x,y\in\Lambda_N$ are such that
$l:=|x-y| \le\min(d(x, \partial\Lambda_N), d(y, \partial\Lambda
_N))$, then one has
%
%
\begin{equation}
\label{eq.2point} \EFK{} {+} {\sigma_x \sigma_y} \le O(1)
\FK{} {+} {0\leftrightarrow\partial\Lambda_{l}}^2 \P^{+} \bigl[\partial(z+ \Lambda_{2l})\leftrightarrow\partial
\Lambda_{N}\bigr],
\end{equation}
where $z$ is the midpoint between $x$ and $y$.
If, on the other hand, one of the points is close to the boundary, in
the sense that
$|x-y| > \min(d(x, \partial\Lambda_N),\break  d(y, \partial\Lambda_N))$,
then one can dominate
$\EFK{}{+}{\sigma_x \sigma_y}$ by $O(1)\FK{}{+}{x\leftrightarrow
\partial\Lambda_{N}} \times\break \FK{}{+}{y\leftrightarrow\partial\Lambda_{N}}$.
\end{pf*}
\end{appendix}

\section*{Acknowledgments}
We wish to thank Mikael De La Salle, Julien Dub\'edat, Hugo
Duminil-Copin, Cl\'ement Hongler and Alain-Sol Sznitman
for useful discussions, an anonymous referee for many useful comments,
and Wouter Kager for Figure~\ref{fig:loops}.
The authors also wish to thank the IHP (Paris) and the NYU Abu Dhabi
Institute for their hospitality.



%

\printaddresses

\end{document}